\documentclass[fleqn,10pt]{wlscirep}

\usepackage{graphicx}
\usepackage[T1]{fontenc}
\usepackage{amsmath,amssymb,amsthm}
\usepackage{physics}
\usepackage{mathtools}
\usepackage{latexsym}
\usepackage[per-mode=symbol]{siunitx}
\usepackage[ruled]{algorithm2e}
\usepackage{here}
\usepackage{float}
\usepackage{lmodern}      % Japanese only
\usepackage{bxpapersize}
\usepackage{subcaption}
\usepackage[version=4]{mhchem}
\usepackage{booktabs}

\usepackage{color}

\newcommand{\argmin}{\mathop{\rm argmin}\limits}

\newcommand{\drv}[2]{\frac{\rmd#1}{\rmd#2}}

\DeclareMathOperator{\diag}{diag}

\DeclareMathOperator*{\minimize}{minimize}

\newcommand{\bbN}{{\mathbb N}}

\newcommand{\bbR}{{\mathbb R}}

\newcommand{\bfb}{{\mathbf b}}

\newcommand{\bfx}{{\mathbf x}}

\newcommand{\bfsigma}{{\pmb{\sigma}}}

\newcommand{\bfA}{{\mathbf A}}

\newcommand{\bfI}{{\mathbf I}}

\newcommand{\bfQ}{{\mathbf Q}}

\newcommand{\bmc}{\mbox{\boldmath $c$}}

\newcommand{\rmd}{{\mathrm d}}

\newcommand{\tilb}{{\tilde b}}

\newcommand{\tilp}{{\tilde p}}
\newcommand{\tilq}{{\tilde q}}

\newcommand{\tilA}{{\tilde A}}

\title{Traffic signal optimization in large-scale urban road networks: an adaptive-predictive controller using Ising models}

\author[1,*,+]{Daisuke Inoue}
\author[2,3,+]{Hiroshi Yamashita}
\author[2]{Kazuyuki Aihara}
\author[1]{Hiroaki Yoshida}
\affil[1]{Toyota Central R\&D Labs., Inc., Nagakute, Aichi 480-1192, Japan}
\affil[2]{International Research Center for Neurointelligence, The University of Tokyo, Bunkyo-ku, Tokyo 113-0033, Japan}
\affil[3]{Present address: Graduate School of Information Science and Technology, Osaka University, 1-5 Yamada-oka, Suita, Osaka 565-0871, Japan.}

\affil[*]{daisuke-inoue@mosk.tytlabs.co.jp}
\affil[+]{These authors contributed equally to this work.}


\begin{abstract}
  Realizing smooth traffic flow is important for achieving carbon neutrality. Adaptive traffic signal control, which considers traffic conditions, has thus attracted attention. However, it is difficult to ensure optimal vehicle flow throughout a large city using existing control methods because of their heavy computational load.
  Here, we propose a control method called AMPIC (Adaptive Model Predictive Ising Controller) that guarantees both scalability and optimality. 
  The proposed method employs model predictive control to solve an optimal control problem at each control interval with explicit consideration of a predictive model of vehicle flow.
  This optimal control problem is transformed into a combinatorial optimization problem with binary variables that is equivalent to the so-called Ising problem.  This transformation allows us to use an Ising solver, which has been widely studied and is expected to have fast and efficient optimization performance.
  We performed numerical experiments using a microscopic traffic simulator for a realistic city road network. The results show that AMPIC enables faster vehicle cruising speed with less waiting time than that achieved by classical control methods, resulting in lower \ce{CO2} emissions. The model predictive approach with a long prediction horizon thus effectively improves control performance.
  Systematic parametric studies on model cities indicate that the proposed method realizes smoother traffic flows for large city road networks. 
  Among Ising solvers, D-Wave's quantum annealing is shown to find near-optimal solutions at a reasonable computational cost.

\end{abstract}

\begin{document}

\flushbottom
\maketitle
% * <john.hammersley@gmail.com> 2015-02-09T12:07:31.197Z:
%
%  Click the title above to edit the author information and abstract
%
\thispagestyle{empty}

% \noindent Please note: Abbreviations should be introduced at the first mention in the main text - no abbreviations lists. Suggested structure of main text (not enforced) is provided below.

\section{Introduction}

\begin{figure}[ht]
  \centering
  \includegraphics[width=\textwidth]{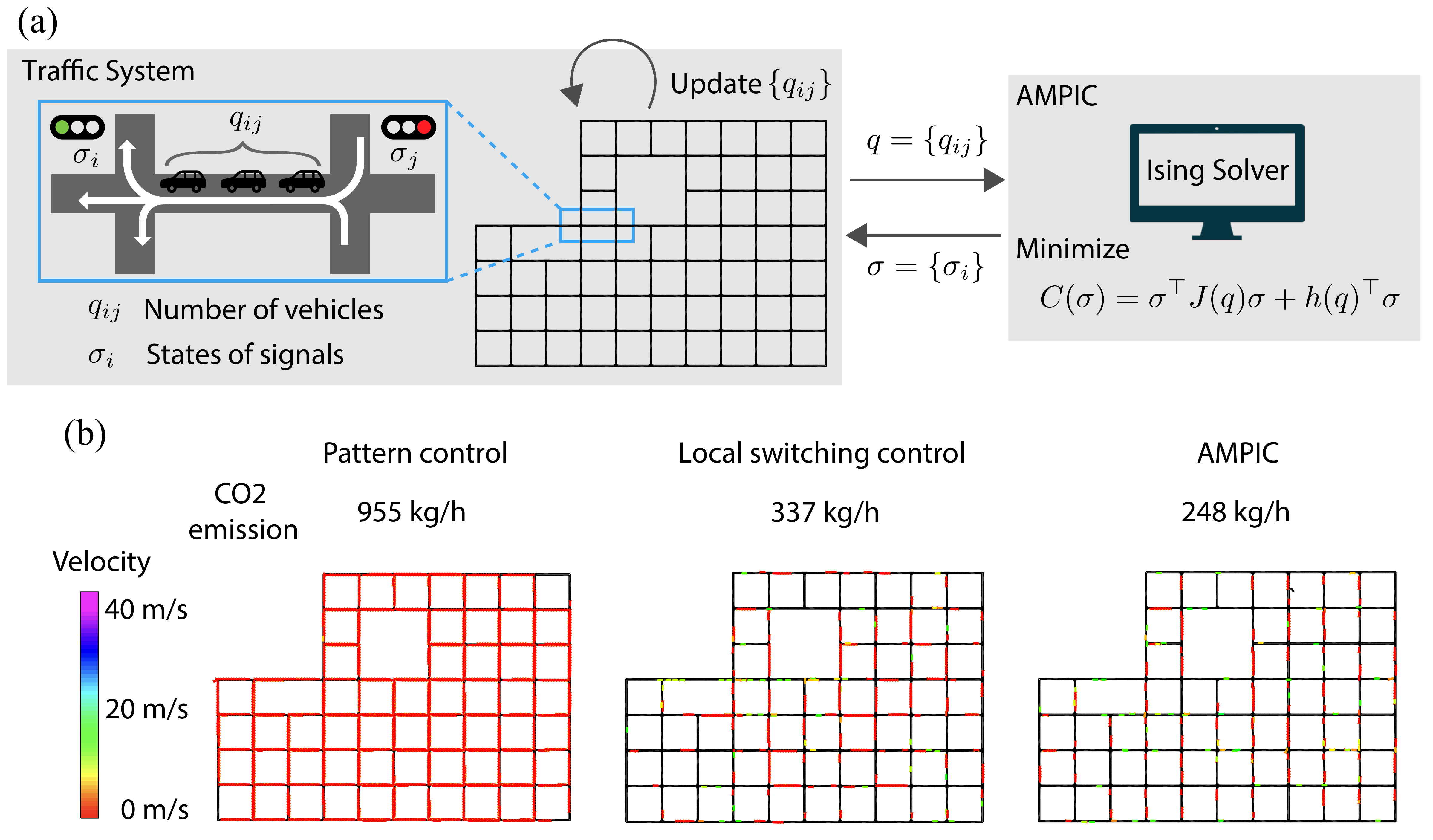}
  \caption{Schematic diagram of the proposed method (AMPIC) and snapshots of vehicles in SUMO with various traffic control methods.
    (a) Vehicle information collected from the traffic system is sent to AMPIC and control quantities computed by AMPIC are sent back (see Section \ref{sec:method-true} for explanation of variables).
    % Vehicle information computed by SUMO is sent to AMPIC and control quantities computed by AMPIC are sent to SUMO (see Section \ref{sec:method-true} for explanation of variables).
    (b) Snapshots of vehicles in numerical evaluation with SUMO on roads that replicate urban areas in Sapporo, Hokkaido, Japan.
    % Snapshots of vehicles on roads that replicate urban areas in Sapporo, Hokkaido, Japan.
    The vehicle generation rate was set to $2.22$ vehicles per second. The snapshots represent the state of vehicles one hour after the start of the simulation.
    The color of a vehicle represents its speed.
    The results for traditional pattern control, local switching control, and control with AMPIC are shown.
    The total \ce{CO2} emissions for each control method are also shown.
  }
  \label{fig:graphic}
\end{figure}

With global economic growth, the amount of traffic in urban areas has continued to increase.
Addressing traffic congestion is essential to mitigate economic losses~\cite{Weisbrod2003Measuring}, ensure driver well-being~\cite{Levy2010Evaluation}, and reduce carbon emissions~\cite{Barth2008RealWorld}.
In particular, the serious impact of \ce{CO2} emissions on global temperatures has been confirmed by climate models and measurement data, requiring immediate action~\cite{Arora2013Carbon,Allen2009Warming,RN2}.
Efficient management and operation of traffic signals are essential for facilitating smooth traffic flow~\cite{Qadri2020Stateofart,Wei2020Survey,Papageorgiou2003Review}.
Signal parameters such as cycle length, green time, and the change interval have traditionally been determined using statistical information based on observed traffic data~\cite{Miller1963Settings,Gartner1983OPAC}.
Recently developed intelligent traffic systems allow the acquisition of real-time traffic information~\cite{Zhao2012Computational,Wu2009development}, enabling adaptive control that dynamically determines signal states based on current traffic conditions~\cite{Cheng2015D2D,Zhang2011DataDriven,Faouzi2011Data,Mirchandani2001realtime,Wang2018Review,Jing2017Adaptive,Lin2011Fast}.
% Recently developed intelligent traffic systems provide real-time traffic information that is useful for determining the optimal state of traffic signals~\cite{Zhao2012Computational,Wu2009development}. 
% Adaptive control technology allows traffic systems to respond dynamically to the current traffic conditions~\cite{Cheng2015D2D,Zhang2011DataDriven,Faouzi2011Data,Mirchandani2001realtime,Wang2018Review,Jing2017Adaptive,Lin2011Fast}.
Several approaches have been proposed to realize adaptive control, including heuristic optimization methods such as genetic algorithms~\cite{Li2015Hybrid}, evolutionary computation~\cite{Chuo2017Evolvable}, and metaheuristic optimization~\cite{Araghi2017Influence,Simoni2019semianalytic}, as well as artificial intelligence models such as neural networks and reinforcement learning~\cite{Rasheed2020Deep, Greguric2020Application, Nishi2018Traffic}.
These adaptive methods are particularly effective for urban traffic networks with rapidly changing traffic patterns~\cite{Wang2018Review}.
However, the computational complexity of determining many control variables and model parameters altogether restricts its use to cases with a network of only a few intersections.
Decentralized control approaches, which split the traffic model and the signal controller, have been proposed~\cite{vandeWeg2019Hierarchical,Mohebifard2019Optimal,Varaiya2013Max}, but the information available for control is limited to the local neighborhood, preventing global optimization throughout the entire city.
Therefore, developing an algorithm that simultaneously and optimally determines control variables within a large city with many intersections is an important challenge~\cite{Qadri2020Stateofart,Choi2016Field}.
% However, the problem of determining control variables is classified as a combinatorial optimization problem, and its computational complexity restricts its use to cases where only a few intersections are handled simultaneously.
% Therefore, extending these algorithms to large cities with a very large number of intersections is an important challenge~\cite{Qadri2020Stateofart,Choi2016Field}.
% Decentralized control approaches, which split the traffic model and the signal controller, have been proposed to address this problem~\cite{vandeWeg2019Hierarchical,Diakaki2002multivariable,Mohebifard2019Optimal,Varaiya2013Max}. 
% However, the information available for control is limited to the local neighborhood, preventing global optimization throughout the entire city.

Recently, methods using Ising models have been actively studied as efficient approaches to solving large-scale combinatorial optimization problems.
The Ising model is a statistical physics model of ferromagnetism that describes the relation between the microscopic state of a spin system and the macroscopic phenomena of magnetic phase transitions~\cite{Yang1952Spontaneous,McCoy2014twodimensional}.
This model is known to be equivalently convertible to a combinatorial optimization problem with binary variables, which includes various practical engineering problems such as delivery planning, packing, and the traveling salesman problem~\cite{Glover2022Applications,Glover2019Tutorial}.
Specialized solvers for finding the ground state of the Ising model, called Ising solvers, have been developed~\cite{Johnson2011Quantuma,Goto2019Combinatorial,Hamerly2019Experimental}.
These solvers exploit the phenomenon of physical quantities reaching a stationary point in the search for a solution and are thus expected to have higher performance in solution seeking than Neumann-type computers.
In particular, a method called quantum annealing, which uses quantum fluctuations to search for the optimal solution, has attracted much attention since D-Wave Systems Inc. provided a commercially available solver, and its performance for various engineering problems has been evaluated~\cite{Neukart2017Traffic,Stollenwerk2020Quantum,Inoue2020Model,Otaki2023Experimental,Okada2023Design,Terada2018Ising,Ohzeki2019Control,Tabi2021Evaluation,Yarkoni2021Multicar}.

Contrary to the expectation suggested by these studies, the application of this computational technology also has a drawback: we need to express the problem in a special form called quadratic unconstrained binary optimization (QUBO). This requires modeling the problem with binary decision variables and an objective function with linear and quadratic terms, and it was a challenge for traffic control problems. In our previous study~\cite{Inoue2021Traffic}, the QUBO model of this problem is proposed by focusing on the balance of the number of vehicles between the roads connected to the same intersection. In addition, the performance of applying Ising solvers to this problem is evaluated with a network of up to thousands of intersections, and the conceptual effectiveness of the method to large-scale networks is demonstrated. In this evaluation, the collective behavior appearing as a spatial pattern is also observed in the optimized traffic lights, suggesting the necessity of global and simultaneous consideration of many traffic lights for optimal control. However, the network and traffic model considered in the study were highly simplified and the applicability to more realistic road networks is not fully investigated.

% This study proposes a global adaptive control algorithm that simultaneously determines many signals in a large city and adapts to the observed traffic conditions.
% To ensure the scalability of the algorithm, we reformulate the optimization problem as a combinatorial optimization problem with binary variables, which is equivalent to the so-called Ising problem (i.e., the problem of finding the ground state of the Ising model).
% The Ising model is a statistical ferromagnetic physics model that describes the relation between the microscopic state of a spin system and the macroscopic phenomena of magnetic phase transitions~\cite{Yang1952Spontaneous,McCoy2014twodimensional}.
% Specialized solvers for finding the ground state of the Ising model, called Ising solvers, have been developed~\cite{Johnson2011Quantuma,Goto2019Combinatorial,Hamerly2019Experimental}.
% These solvers are expected to provide fast solutions even for large-scale optimization problems.
% 
% 
In this paper, we propose a method to achieve global adaptive control that simultaneously determines many signals in a large city and more flexibly adapts to road network and observed traffic conditions by solving the Ising problem.
The main contribution of this study is summarized as follows:
% The main features of the method proposed in the present study, called AMPIC (Adaptive Model Predictive Ising Controller), are summarized as follows:
\begin{itemize}
  \item
        The proposed control method uses an Ising solver to determine the traffic signal state at each time.
        % The proposed control method transforms the problem of determining traffic signals at each time into an Ising problem, and the search for optimal variables is carried out by the Ising solver. 
        The optimization problem of traffic signals using an internal model that predicts the flow dynamics of vehicles traveling on an arbitrary road network (up to four-way intersections) is shown to be equivalent to the Ising problem.
        This allows scalable traffic signal control even in large-scale road networks.

  \item
        The proposed method works adaptively according to traffic conditions.
        The parameters of the internal model used by the proposed control method are adaptively updated using observed vehicle flow information.
        This approach reduces modeling error (compared to actual traffic flow) and allows the controller flexibility to adapt to various traffic scenarios.
        % The proposed method works adaptively according to traffic conditions. 
        % It includes an internal model that predicts the flow dynamics of vehicles traveling in an arbitrary road network (with up to four-way intersections). 
        % The parameters of the internal model are updated adaptively using vehicle flow information observed during each control cycle. 
        % This approach reduces modeling errors (compared to actual traffic flow) and allows the controller to respond flexibly to various traffic scenarios.

  \item Model predictive control is employed. The controller's internal model is used to predict traffic conditions up to multiple control cycles ahead and minimize an objective function to improve future traffic conditions.
        This approach avoids short-sighted control and is expected to achieve advanced traffic management.
\end{itemize}
We name the proposed method AMPIC (Adaptive Model Predictive Ising Controller) after the above-described features.

% Previous studies on signal control problems have related the optimization problem with Ising-type problems
Including our previous study~\cite{Inoue2021Traffic}, there are several previous studies that relate signal control problems to Ising-type problems.~\cite{Hussain2020Optimal,Marchesin2023Improving,Shikanai2023Traffic}.
% Our previous paper~\cite{Inoue2021Traffic} proposed a method for predicting and optimizing the traffic conditions one control cycle ahead for a group of vehicles traveling on a road with lattice periodic boundary conditions.
% and verified the performance using the same internal model of the controller.
Hussain et al.~\cite{Hussain2020Optimal} proposed a method to maximize traffic flow in a lattice network by coordinating adjacent signals.
% and verified the performance using vehicle models that travel according to the given signal indications.
Marchesin et al.~\cite{Marchesin2023Improving} considered an optimization problem and proposed a method whose parameters are adjusted based on vehicle flow observations for an arbitrary group of vehicles and a road network.
% They verified the performance of their method using the same model as the controller's internal model.
Shikanai et al.~\cite{Shikanai2023Traffic} considered an optimization problem and proposed a method to suppress traffic congestion that adapts to the states of vehicles of arbitrary vehicles and road networks.
% They verified the performance of their method using SUMO.
%The method proposed by Shikanai et al.~\cite{Shikanai2023Traffic} is similar to that proposed in the present study.
Except for Ref.~\citeonline{Shikanai2023Traffic}, the validation of controllers is performed on coarse-grained macroscopic models, while ours is based on more realistic microscopic simulations in which individual vehicles move according to traffic rules.
In contrast to Ref.~\citeonline{Shikanai2023Traffic},
% our objective function is designed such that the optimization problem becomes compatible with the Ising model.
% , which eliminates ambiguous hyper-parameters. 
our method employs model predictive control, which enables more advanced signal control by predicting future traffic conditions.
Several methods have already been proposed to use model predictive control for traffic signal optimization~\cite{Lin2011Fast,Ye2019survey,Mirchandani2001RealTimea,Guo2019Urban}.
While these methods typically require a large amount of computation as the network size increases, our approach through the Ising model alleviates this difficulty.
% Our objective function, however, is designed such that the optimization problem becomes compatible with the Ising model, which eliminates ambiguous hyper-parameters. In addition, our method applies model predictive control.

In the present study,
the performance of AMPIC is systematically evaluated using an external microscopic traffic flow simulator, namely Simulation of Urban MObility (SUMO)~\cite{SUMO2018}, which is widely used to model realistic urban traffic (see Fig.~\ref{fig:graphic}).
The performance is evaluated in terms of practical measures such as waiting time and \ce{CO2} emissions.
% AMPIC is proposed and applied to a realistic urban road network.
% Numerical experiments are performed using SUMO. 
The results show that AMPIC increases the vehicle's cruising speed and reduces waiting times compared with those obtained with conventional control methods. They also show that it considerably reduces \ce{CO2} emissions.
Notably, the model predictive approach significantly enhances control performance, especially when it employs an extended prediction horizon.
Parameter studies with model cities indicate that AMPIC achieves more efficient traffic flows compared to those obtained using other control protocols, particularly for very large cities. We also show that D-Wave's quantum annealing consistently finds near-optimal solutions at reasonable computational cost.

\section{Method}\label{sec:method-true}

\subsection{Model predictive control of traffic signals}\label{sec:problem_formulation}

In this study, we consider the problem of controlling traffic signals in a road network with $N$ intersections. We consider the graph $G=(V, E)$, where $V=\{1,\ldots,N\}$ denotes the index of an intersection and $E$ denotes the roads connecting to it.
In this study, we assume that the network is directed; a road from intersection $j$ to $i$ is denoted by $(i,j)$.
% If there is an edge $(i,j)$, we assume that there is also an opposite edge $(j,i)$.
We consider up to 4-way intersections (i.e., there are at most four roads leading to intersection $i$).

\begin{figure}[ht]
  \centering
  \includegraphics[width=0.8\textwidth]{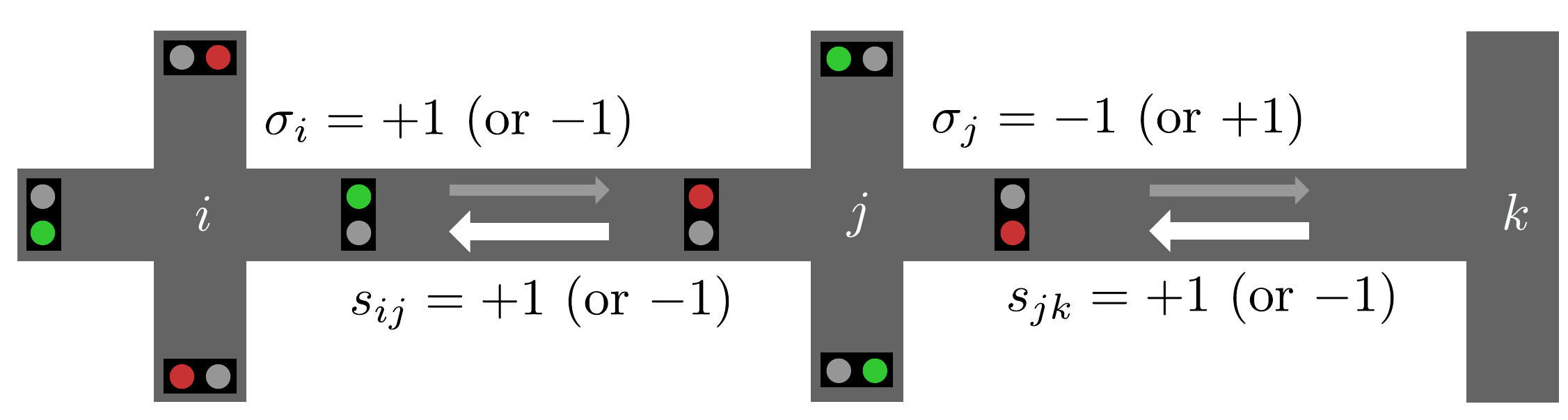}
  \caption{Definition of variables for the road network.}
  \label{fig:sij}
\end{figure}

Let $\sigma_i(t)$ denote the state of the traffic signal at intersection $i$ at time $t$.
We assume that each traffic signal has two states, namely $\sigma_i\in\{\pm1\}$, and that it can be either red or green depending on its state.
% We use $s_{ij}$ to represent the state of the traffic light (red or green) at intersection $i$ for road $(i, j)$. 
Based on traffic light states and road positioning, we introduce the variable $s_{ij}\in\{\pm 1\}$, that is, if the traffic light $i$ with the state $\sigma_i$ indicates green for the road $(i, j)$, we assign the state value $\sigma_i$ to $s_{ij}$.
Then, the traffic light indicates green for the road $(i,j)$ when $\sigma_i s_{ij}=+1$ and red when $\sigma_i s_{ij}=-1$.
The definitions of these symbols are shown in Fig.~\ref{fig:sij}.
We assume that the signal state $\sigma_i(t)$ is determined at discrete time $t=\tau k\ (k\in \bbN)$ with a predetermined control cycle $\tau>0$.
Once the signal is determined, the state $\sigma_i(t)$ is fixed until the next control cycle $\sigma_i(t+\tau)$.
The number of vehicles on road $(i,j)$ at time $t$ is denoted by $q_{ij}(t)$.

Depending on the indication of the traffic light, the roads leading to the intersection are classified as being on the red side or the green side.
We define vehicle bias $x_i$ at intersection $i$ to represent how the situation differs from the condition where vehicles entering the intersection are equally distributed on the red and green sides:
\begin{align}\label{eq:def_x}
  % x_i(t) = 2\sum_{j\in N_i} \eta_{ij}s_{ij}q_{ij}(t),
  x_i(t) = \sum_{j\in N_i} \eta_{ij}s_{ij}q_{ij}(t),
\end{align}
where $N_i$ is the set of intersection indices that have a road leading to intersection $i$.
The parameter $\eta_{ij}>0$ is a weight parameter determined based on the road and intersection geometry. The setting of this parameter is described in the supplementary materials.
In a typical example of east-west and north-south road intersections, $x_i$ represents the difference between the number of vehicles present in the east-west direction and that in the north-south direction.

The vectors $x(t)\coloneqq [x_1(t), x_2(t), \ldots, x_N(t)]^\top\in\bbR^N$ and $\sigma(t)\coloneqq [\sigma_1(t), \sigma_2(t), \ldots, \sigma_N(t)]^\top\in\bbR^N$ are defined as the vectors of the vehicle bias and traffic signal state, respectively.
Based on bias vector $x(t)$ at a certain time $t$, we find the present $\sigma(t)$ by minimizing the following objective function:
\begin{align}\label{eq:eval_func_multi}
  \begin{split}
    C(\sigma([t,\ldots,t+\tau k_{\mathrm h}])) & = \sum_{k=1}^{k_{\mathrm h}} x(t+\tau k)^\top Qx(t+\tau k),
    % - \sum_{k=0}^{k_{\mathrm h}-1}\qty{ w_k \sigma(t-1+k)^\top R \sigma(t+k)}
  \end{split}
\end{align}
where $k_{\mathrm h}\in\bbN$ is a constant called the prediction horizon; a large value of $k_{\mathrm h}$ takes into account a long-term future. Diagonal matrix $Q\in\bbR^{N\times N}$ is used to set the relative weight of importance on the intersections; $Q=I$, where $I$ is the identity
matrix, means that all intersections are equally important.
Equation~\eqref{eq:eval_func_multi} aggregates the vehicle bias for every intersection at each time step.
By minimizing this quantity, we expect to eliminate the uneven flow of vehicles throughout the city, resulting in smooth, congestion-free traffic flow.

Equation~\eqref{eq:eval_func_multi} includes future information after the current time $t$.
To predict this information, we need to model the time evolution of vehicle bias $x(t)$.
In this study, we assume that the changes in the number of vehicles $q(t)$ are constant given the traffic light pattern $\sigma$ within the prediction horizon, and the bias can be expressed by the following linear difference equation:
\begin{align}\label{eq:dynamics_x_vec_in_result}
  x(t+\tau) & = x(t) + \tilde{A}{\sigma}(t) + \tilde{b},
\end{align}
where the matrix $\tilde A$ and vector $\tilde b$ are determined by the statistical information collected from the road network.
The derivation of Eq.~\eqref{eq:dynamics_x_vec_in_result} and the specific forms of $\tilde A$ and $\tilde b$ are described in the supplementary materials.
By substituting the predictive model \eqref{eq:dynamics_x_vec_in_result} into the objective function \eqref{eq:eval_func_multi}, we obtain the so-called Ising problem with $\sigma(s)\ (s\in\{t,\ldots,t+\tau k_{\mathrm h}\})$ as the decision variable, which is discussed in the following subsection.

\subsection{Ising models, problems, and solvers}

By combining Eqs.~~\eqref{eq:eval_func_multi} and \eqref{eq:dynamics_x_vec_in_result}, we can rewrite the above optimal control problem in the following form:
\begin{align}
  \begin{aligned}\label{eq:QUBO}
    \minimize_\sigma & \;\sum_{i<j}J_{ij}\sigma_i\sigma_j + \sum_{i} h_i\sigma_i,
    % \\
    % \text{s.t}&\;\sigma_i\in\{-1,1\}\;\forall\;i,
  \end{aligned}
\end{align}
where we consider $N k_{\mathrm h}$ decision variables, denoted by $\sigma=[\sigma_1, \ldots, \sigma_{N k_{\mathrm h}}]^\top\in\{-1,1\}^{N k_{\mathrm h}}$, and consider $J\in\bbR^{{N k_{\mathrm h}}\times {N k_{\mathrm h}}}$ and $h\in\bbR^{N k_{\mathrm h}}$ as the parameters that characterize this problem.
The specific definitions of $J$ and $h$ in terms of the model parameters such as $\tilde A$, $\tilde b$, and $Q$ are given in the supplemental materials.

In this problem, the decision variables are binary, either $-1$ or $+1$, and the objective function is quadratic.
Since these variables and objective function correspond to the spin variables and the energy of the Ising model in statistical physics, respectively,
% the problem in Eq.~\eqref{eq:QUBO} is regarded as a problem of finding a spin state that minimizes the energy of the Ising model.
% Because of this correspondence, 
we consider our problem to be an Ising problem.

The Ising problem is both an engineering optimization problem and a physical model.
Researchers have thus attempted to utilize the model as a solver for optimization problems by constructing a physical system or hardware that minimizes Eq.~\eqref{eq:QUBO} (see, for example, Refs.~\citeonline{Johnson2011Quantuma,Goto2019Combinatorial,Hamerly2019Experimental}).
In this study, these solvers are referred to as Ising solvers.

Examples of such hardware include coherent Ising machines~\cite{Inagaki2016coherent,Hamerly2019Experimental}, simulated bifurcation machines~\cite{Goto2019Combinatorial}, digital annealers~\cite{Matsubara2018IsingModel,Aramon2019PhysicsInspired}, and quantum annealing machines from D-Wave Systems Inc.\cite{Kadowaki1998Quantum,Johnson2011Quantuma}
Among these, quantum annealing machines have attracted attention as a non-von Neumann, commercially available computer architecture that takes advantage of quantum fluctuations. Their applications are currently under investigation~\cite{King2015Benchmarking,McGeoch2013Experimental,Venturelli2016Quantum,OMalley2018Nonnegative,Ohzeki2018Optimization,Inoue2020Model,Ayanzadeh2020Reinforcementa}.
% Therefore, in this study, we evaluate the performance of the proposed method using multiple Ising solvers including the quantum annealing machine. 
% implemented on a quantum annealing machine. 
We note that all these solvers are capable of performing our proposed method.
% is feasible by can be implemented in all these solvers. 
% The proposed method is compatible with any above solvers.

\subsection{Coupled simulation with SUMO}

We evaluate the proposed control algorithm by applying it to traffic lights that are obeyed by individual vehicles traveling in a city road network.
For evaluation, a microscopic traffic simulator, in which each vehicle travels through the city under realistic traffic conditions, is coupled to the proposed controller.
In particular, we use the open-source software SUMO~\cite{SUMO2018} as the simulator.

As shown in Fig.~\ref{fig:graphic}, SUMO and AMPIC are coupled via variables that are used in their respective calculations.
In SUMO, vehicles travel along a preplanned route, whose origin and destination are randomly chosen.
The position, speed, and other states of the vehicles are updated sequentially every second.
The number of vehicles on roads $q_{ij}$ is calculated and sent to AMPIC.
In AMPIC, the Ising model is constructed using the values of $q_{ij}$ and then the optimal states of all traffic signals are calculated.
SUMO receives the signal states, which are used to update the position of the vehicles.
Note that the calculation of control inputs does not necessarily require simulation with SUMO, but requires only the number of vehicles observed at each time to identify the linear equation \eqref{eq:dynamics_x_vec_in_result}.
See the supplementary materials for the detailed simulation procedures.

\section{Results}

\subsection{Evaluation criteria}
The performance of the proposed controller is compared with that of existing methods using the following performance indicators:
\begin{itemize}
  \item Mean velocity: the average speed of all vehicles in the road network.
        %    Precisely, the average vehicle speed per road is averaged for all roads.
  \item Waiting vehicle ratio: the ratio of the number of stopped vehicles to the number of all vehicles in the road network.
        Here, vehicles moving at a speed below $0.1\ \si{\meter\per\second}$ are counted as stopped vehicles.
  \item \ce{CO2} emissions: the total amount of carbon dioxide emitted by all vehicles on the road per time step, which is computed using a model included in SUMO~\cite{krajzewicz2015second}.
  \item Sum of squared vehicle bias: the value of the objective function defined in Eq.~\eqref{eq:eval_func_multi}, evaluated at each time step using the obtained $\sigma(t)$:
        \begin{align}
          \tilde C(t) = x(t)^\top Q x(t).
        \end{align}
\end{itemize}
Note that the first three performance indicators directly evaluate the desirability of traffic conditions and the last performance indicator is closely related to the minimized objective function (i.e., the Ising model).
We compare the proposed method with the following control methods:
\begin{itemize}
  \item Local switching control~\cite{Suzuki2013Chaotic}:
        This controller greedily reduces the vehicle bias for each intersection, considering only the local vehicle bias rather than minimizing a global objective function.
        The state of signal $\sigma_i$ at each intersection $i$ is defined as
        \begin{align}
          \begin{cases}\label{eq:local_control}
            % \sigma_i(t) \leftarrow +1 & \text{if } x_i(t) \ge +\theta,\\
            % \sigma_i(t) \leftarrow -1 & \text{if } x_i(t) \le -\theta.
            \sigma_i(t) \leftarrow +1                & \text{if } x_i(t) > 0, \\
            \sigma_i(t) \leftarrow \sigma_i(t-\tau ) & \text{if } x_i(t) = 0, \\
            \sigma_i(t) \leftarrow -1                & \text{if } x_i(t) < 0.
          \end{cases}
        \end{align}
  \item Random control:
        This control method switches the state of each signal with a certain probability at each control cycle.
        The switching probability per control cycle is assumed to be $0.5$.
  \item Pattern control:
        This control method switches the state of the signal at each predetermined control cycle.
        The switching cycle is assumed to be once every two control cycles for compatibility with random control.
\end{itemize}

The vehicle origin location, destination, and routing follow SUMO's default settings (see the supplementary materials for details).
Unless otherwise noted, the parameters of the controller are defined as follows: the control period is $\tau = 60\ \si{\second}$, the prediction horizon is $k_{\mathrm h}=1$, and the weight matrix is $Q\equiv I$.
For the Ising solver, we use \emph{SimulatedAnnealingSampler}~\cite{SASampler} provided by D-Wave Systems and the parameter \emph{num\_reads} is set to $1000$.
For the other hyperparameters in the solver, the default values are used.

\begin{figure}[th]
  \centering
  \includegraphics[width=0.7\textwidth]{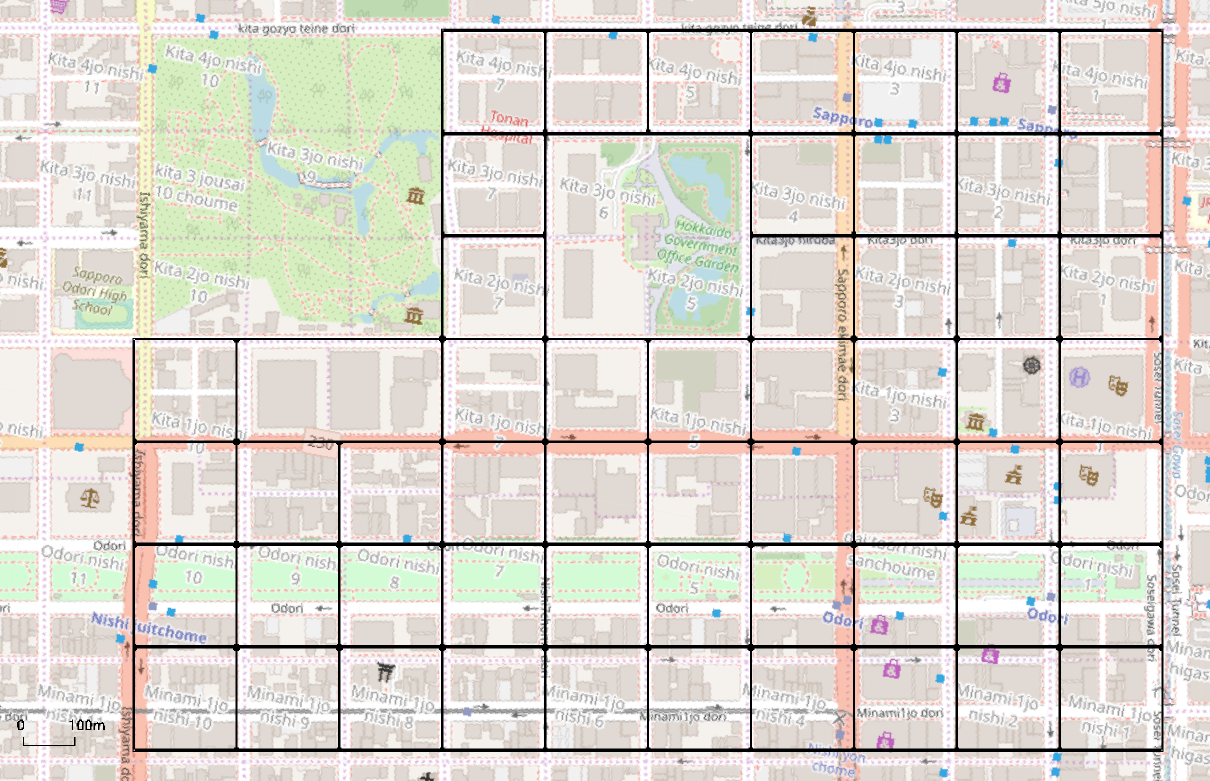}
  \caption{Road network created to replicate an urban area of Sapporo, Hokkaido, Japan.
    %     The map image is taken from OpenStreetMap (\url{https://www.openstreetmap.org}).
    The base map image is provided by \copyright OpenStreetMap~\cite{OSM,OSMPaper,OSMap}.
  }
  \label{fig:sapporo-network}
\end{figure}

\subsection{Results for Sapporo}

We created a road network that corresponds to an urban area in Sapporo, as shown in Fig.~\ref{fig:sapporo-network}.
All intersections have traffic signals. The color of the lights is determined by AMPIC such that the objective function in Eq.~\eqref{eq:eval_func_multi} is minimized at each control cycle.

\begin{figure}[ht]
  \centering
  \begin{subfigure}{0.49\textwidth}
    \includegraphics[width=0.9\linewidth]{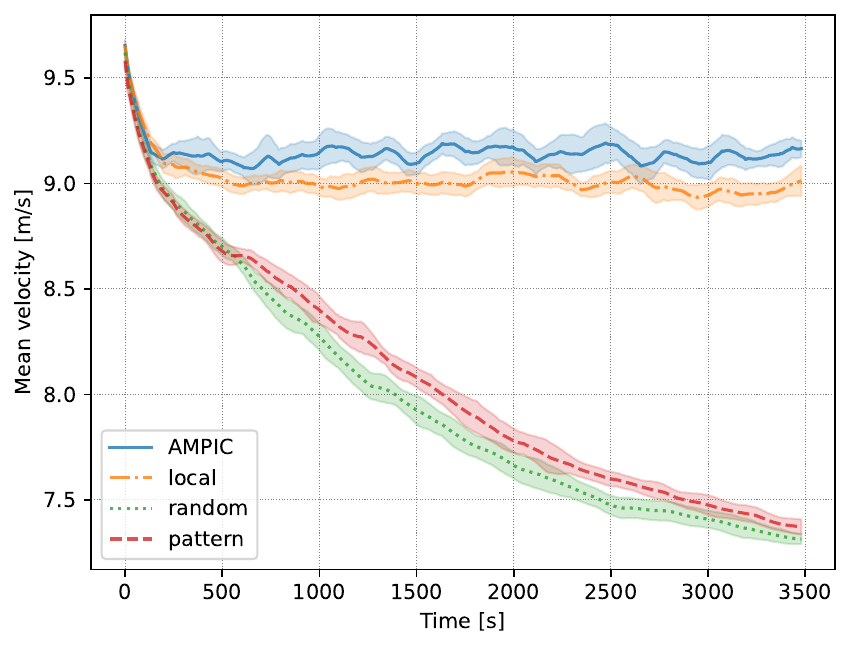}
    \caption{Average velocity of vehicles.}
    \label{fig:time_response_vmean}
  \end{subfigure}
  \begin{subfigure}{0.49\textwidth}
    \includegraphics[width=0.9\linewidth]{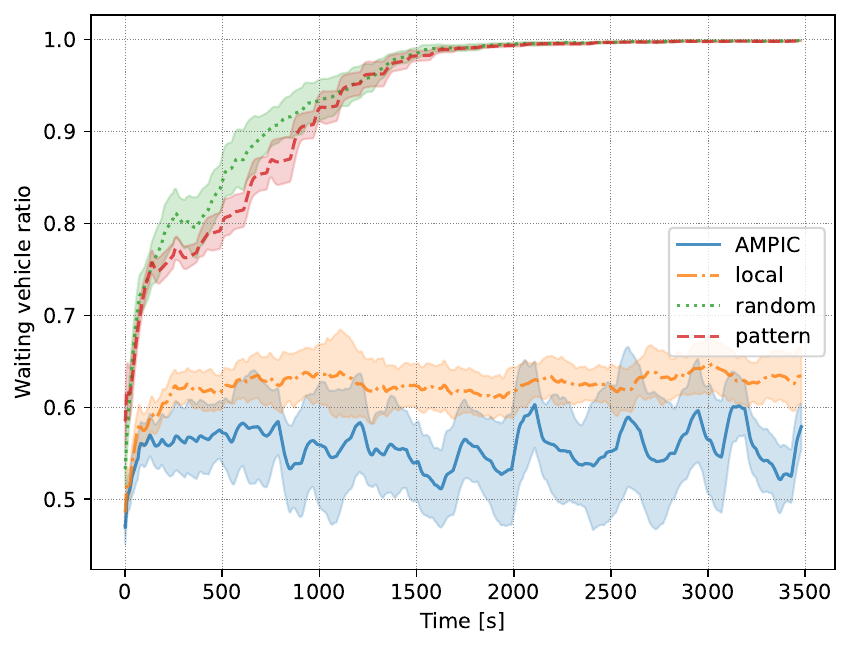}
    \caption{Waiting vehicle ratio.}
    \label{fig:time_response_waitrate}
  \end{subfigure}
  \begin{subfigure}{0.49\textwidth}
    \includegraphics[width=0.9\linewidth]{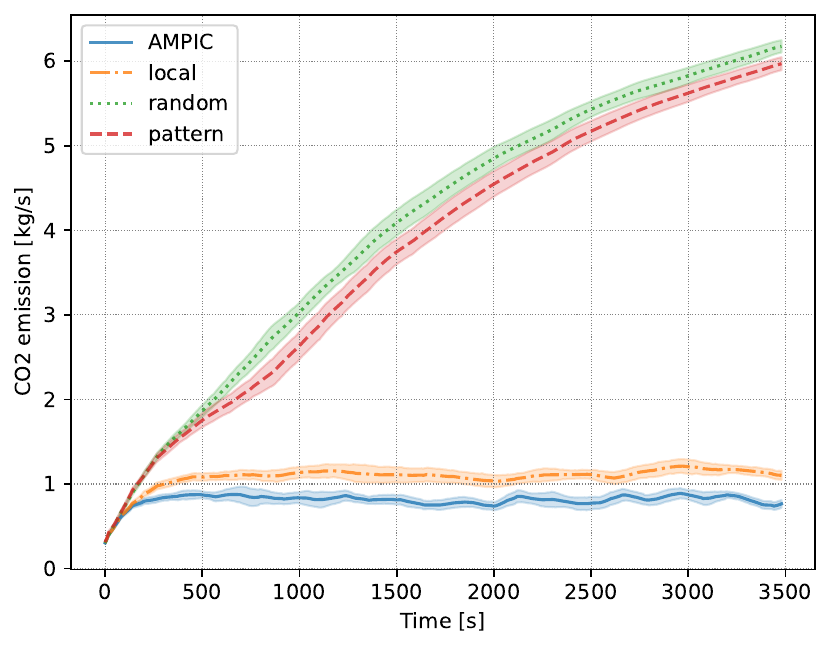}
    \caption{\ce{CO2} emissions.}
    \label{fig:time_response_co2}
  \end{subfigure}
  \begin{subfigure}{0.49\textwidth}
    \includegraphics[width=0.9\textwidth]{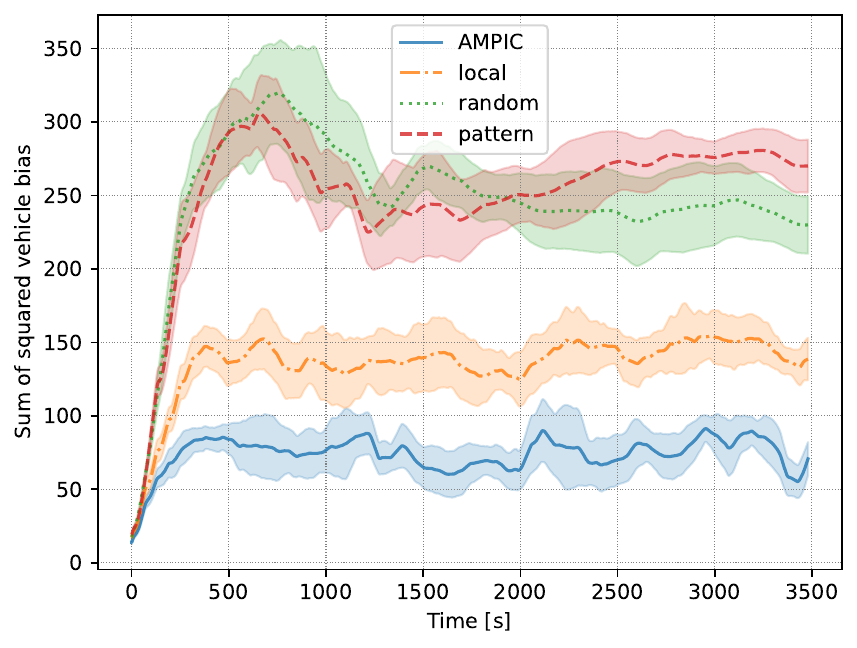}
    \caption{Sum of squared vehicle bias.}
    \label{fig:time_response_energy}
  \end{subfigure}
  \caption{Time evolution of performance indicators for various methods. The proposed method (blue solid) is compared with local control (orange dash-dotted), random control (green dotted), and pattern control (red dashed).
    The colored area represents the standard error of each performance indicator for various seeds of the random number used for vehicle route generation.}
  \label{fig:time_response}
\end{figure}

The time evolution of the mean velocity, waiting vehicle ratio, \ce{CO2} emissions, and the sum of squared vehicle bias are shown in Fig.~\ref{fig:time_response}.
The number of vehicles added to the road network per time step, which we refer to as the vehicle generation rate, is set to $2.22$ vehicles per second.
The figure shows the average value over five simulation runs, in which the vehicles had different routes.
For legibility,
the 120-second moving average for each performance indicator is shown.
As shown in Figs.~\ref{fig:time_response_vmean} and \ref{fig:time_response_waitrate}, with the non-adaptive methods (i.e., random and pattern control methods), the mean velocity decreases over time and the waiting vehicle ratio approaches $1$, which means that congestion occurs.
The congestion results in an increasing number of vehicles in the network. Accordingly, the \ce{CO2} emissions increase steadily, as shown in Fig.~\ref{fig:time_response_co2}.
In contrast, with local switching control and the proposed control method, the mean velocity and waiting vehicle ratio approach a steady state after a short initial transient response.
This results in low \ce{CO2} emissions due to the balance between the number of vehicles that are added and the number of vehicles that are removed (i.e., vehicles that reach their destinations).
The proposed control method achieves shorter waiting times and faster cruising speed
than those obtained with local switching control; the \ce{CO2} emissions are accordingly lower.
In Fig.~\ref{fig:time_response_energy}, the proposed method shows the lowest vehicle bias among all methods. This is expected as this performance indicator is directly minimized in the proposed method, but not in the other methods.

\begin{figure}[ht]
  \centering
  \begin{subfigure}{0.49\textwidth}
    \includegraphics[width=0.9\textwidth]{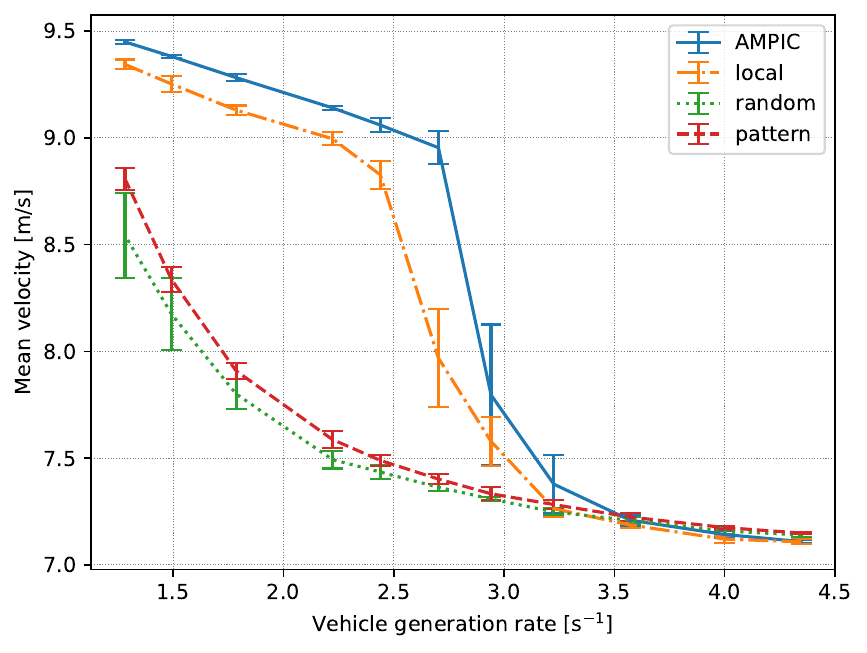}
    \caption{Mean velocity.}
    \label{fig:various_p_speed}
  \end{subfigure}
  \begin{subfigure}{0.49\textwidth}
    \includegraphics[width=0.9\textwidth]{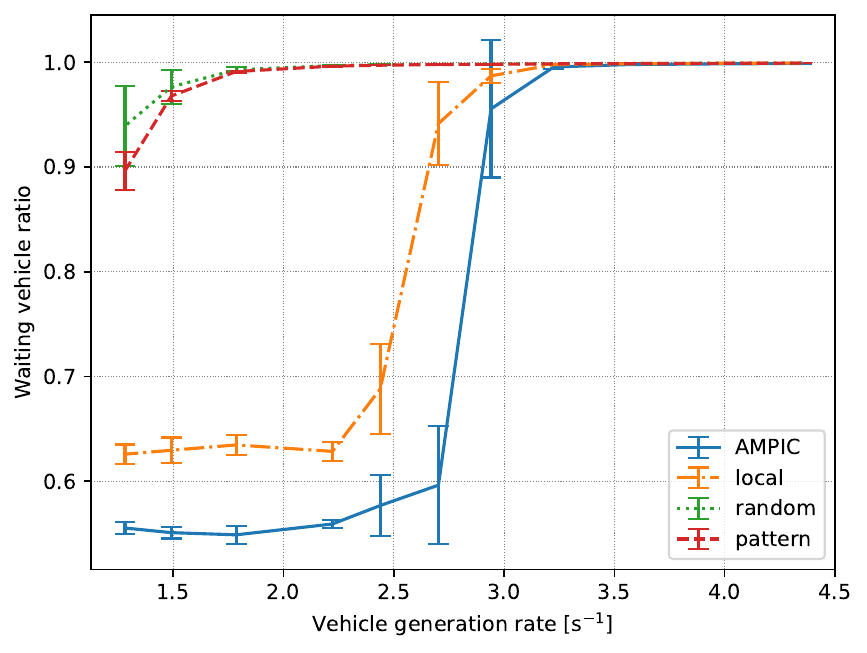}
    \caption{Waiting vehicle ratio.}
    \label{fig:various_p_waitrate}
  \end{subfigure}
  \begin{subfigure}{0.49\textwidth}
    \includegraphics[width=0.9\textwidth]{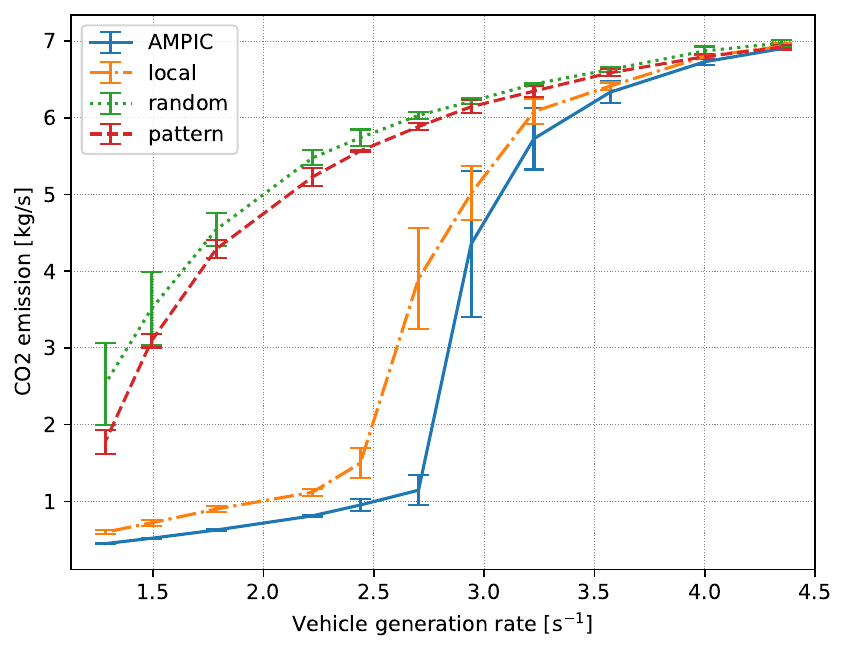}
    \caption{\ce{CO2} emissions.}
    \label{fig:various_p_co2}
  \end{subfigure}
  \begin{subfigure}{0.49\textwidth}
    \includegraphics[width=0.9\textwidth]{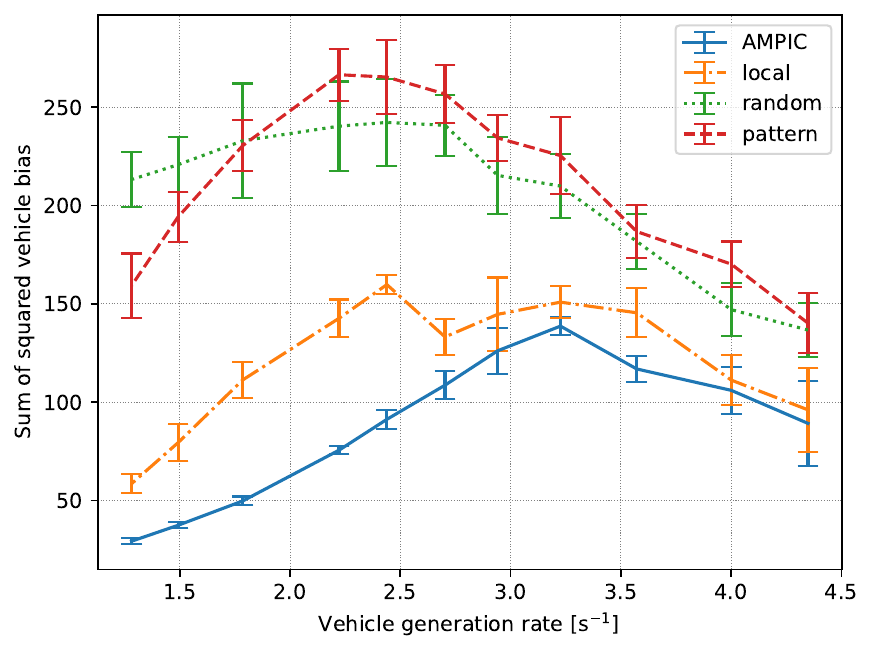}
    \caption{Sum of squared vehicle bias.}
    \label{fig:various_p_energy}
  \end{subfigure}
  \caption{Time-averaged performance indicators for various vehicle generation rates for various methods.
    The proposed method (blue solid) is compared with local control (orange dash-dotted), random control (green dotted), and pattern control (red dashed).
    The error bars represent the standard error of each performance indicator for various seeds of the random number used for vehicle route generation.}
  \label{fig:various_p}
\end{figure}

We next examine situations with different numbers of vehicles. Simulations were performed with various vehicle generation rates.
Figure~\ref{fig:various_p} shows the time-averaged values of the mean velocity, the waiting vehicle ratio, \ce{CO2} emissions, and the sum of squared vehicle bias as functions of the vehicle generation rate.
For each simulation run, the value of a given performance indicator was averaged in time over $3600$ seconds.
Each data point and error bar indicate the mean value and standard error, respectively, for five simulation runs with different vehicle routes.
% The error bars represent the standard error of each statistic.
Figures~\ref{fig:various_p_speed} and \ref{fig:various_p_waitrate}
% Figs.~\ref{fig:various_p} and \ref{fig:various_p_waitrate}
show that in general, an increase in the vehicle generation rate leads to a lower average speed and more stopped vehicles.
The proposed method mitigates the decrease in mean velocity
and thus the waiting vehicle ratio and the \ce{CO2} emissions are suppressed.
Local switching control shows almost the same performance for a large vehicle generation rate. However, note that the considered parameter range corresponds to a highly congested situation, with most vehicles being stopped.
Similar trends were found for \ce{CO2} emissions in Fig.~\ref{fig:various_p_co2}, suggesting that better traffic conditions are also beneficial for decarbonization.
As shown in Fig.~\ref{fig:various_p_energy}, the proposed method minimizes the energy of the Ising model the most for all vehicle generation rates; this result confirms the correspondence between the performance indicator defined via vehicle bias and other performance indicators (velocity, waiting vehicle ratio, and \ce{CO2} emissions).
% that directly reflects the intuitive traffic conditions.

\subsection{Size of road network}

We now examine the performance of the proposed method for road networks with various sizes. Here, we use square lattice networks that consist of $N$ intersections, with $N$ varied from $5\times 5$ to $25\times 25$.
The distance between adjacent intersections is fixed at $100\ \si{\meter}$.
% We focus on mean velocity,  of stopped vehicles, and energy of the Ising model.
Since local switching control was found to be the most comparable method in the previous subsection, it is used here in the comparison with AMPIC.

\begin{figure}[ht]
  \centering
  \begin{subfigure}{0.32\textwidth}
    \centering
    \includegraphics[width=0.9\columnwidth]{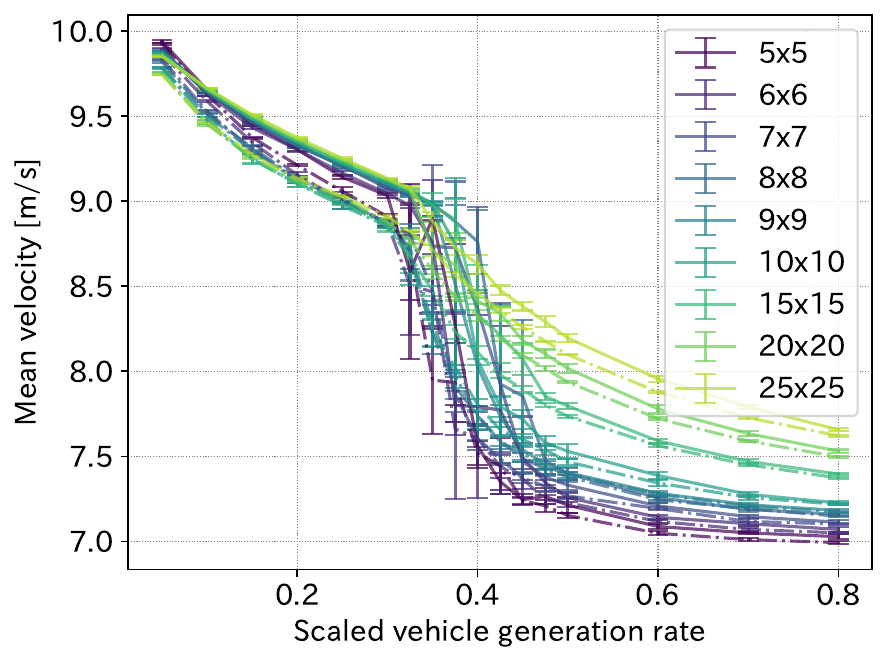}
    \caption{Mean velocity.}
    \label{fig:vmean-networksize}
  \end{subfigure}
  \begin{subfigure}{0.32\textwidth}
    \centering
    \includegraphics[width=0.9\columnwidth]{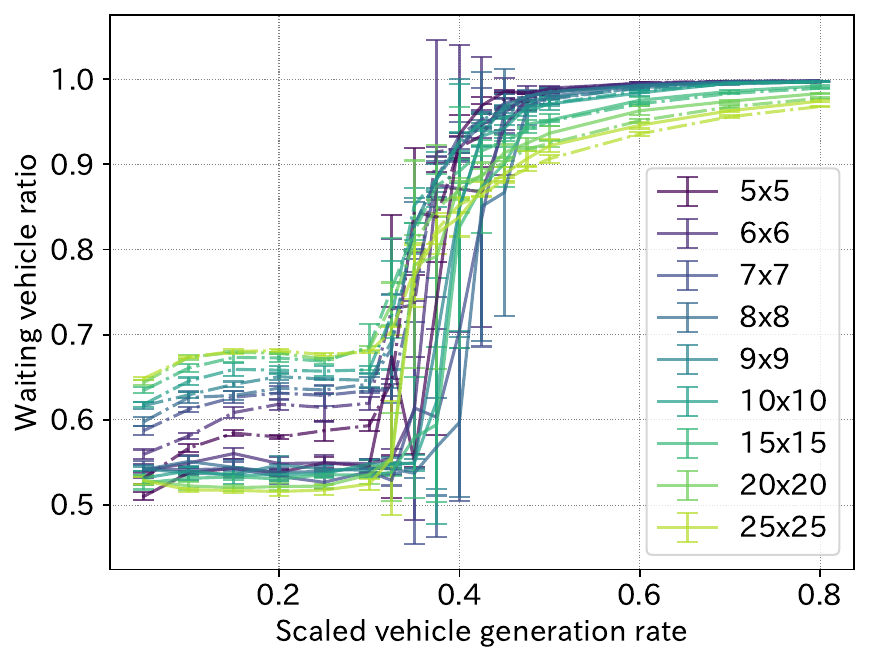}
    \caption{Waiting vehicle ratio.}
    \label{fig:waterate-networksize}
  \end{subfigure}
  \begin{subfigure}{0.32\textwidth}
    \centering
    \includegraphics[width=0.9\columnwidth]{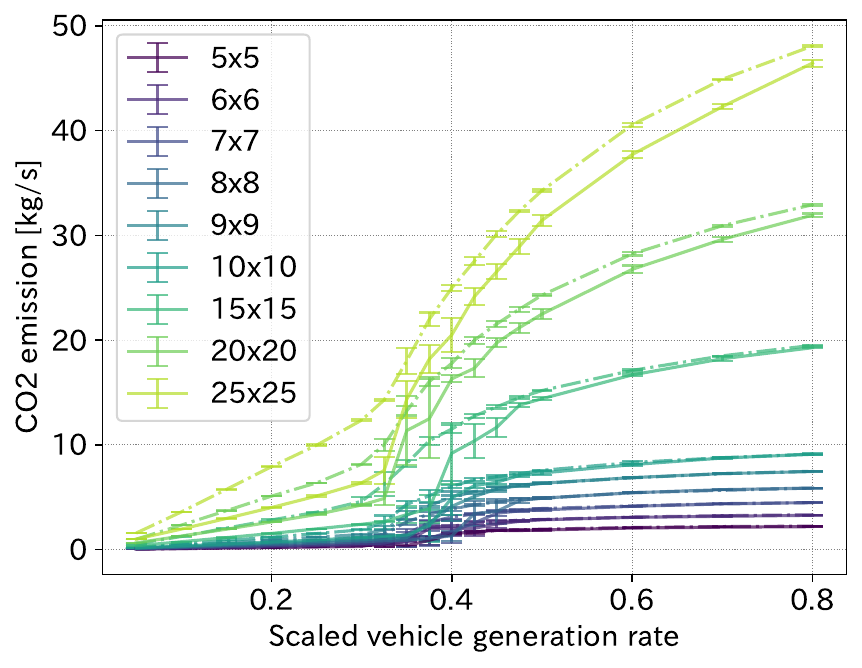}
    \caption{\ce{CO2} emissions.}
    \label{fig:co2-networksize}
  \end{subfigure}
  \begin{subfigure}{0.32\textwidth}
    \centering
    \includegraphics[width=0.9\columnwidth]{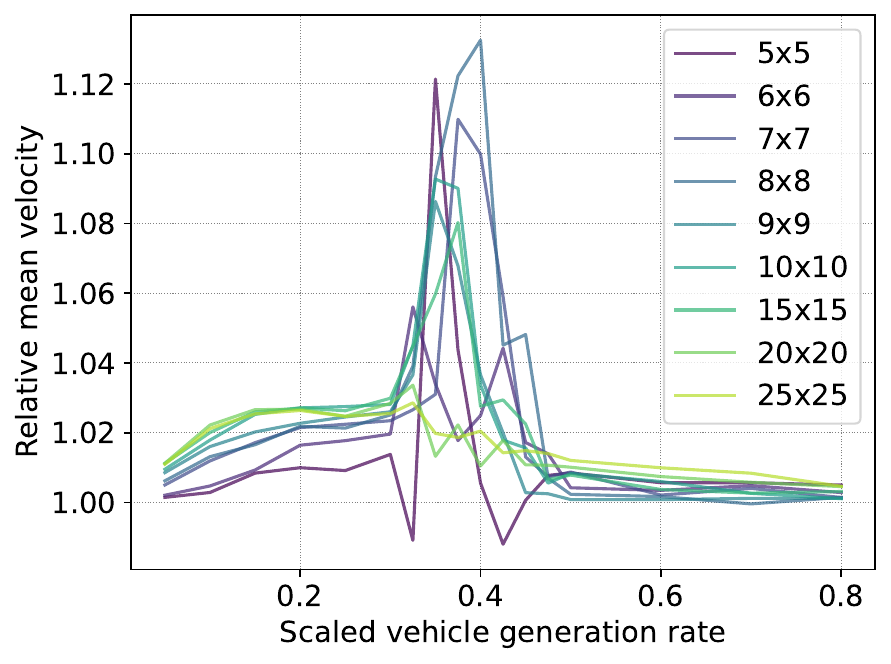}
    \caption{Relative mean velocity.}
    \label{fig:vmean-networksize-rel}
  \end{subfigure}
  \begin{subfigure}{0.32\textwidth}
    \centering
    \includegraphics[width=0.9\columnwidth]{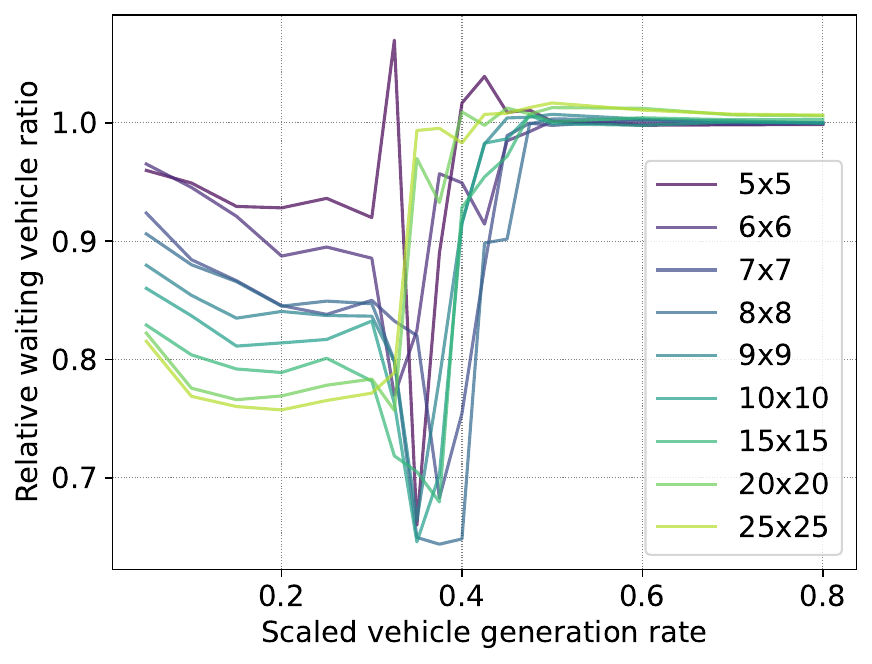}
    \caption{Relative waiting vehicle ratio.}
    \label{fig:waterate-networksize-rel}
  \end{subfigure}
  \begin{subfigure}{0.32\textwidth}
    \centering
    \includegraphics[width=0.9\columnwidth]{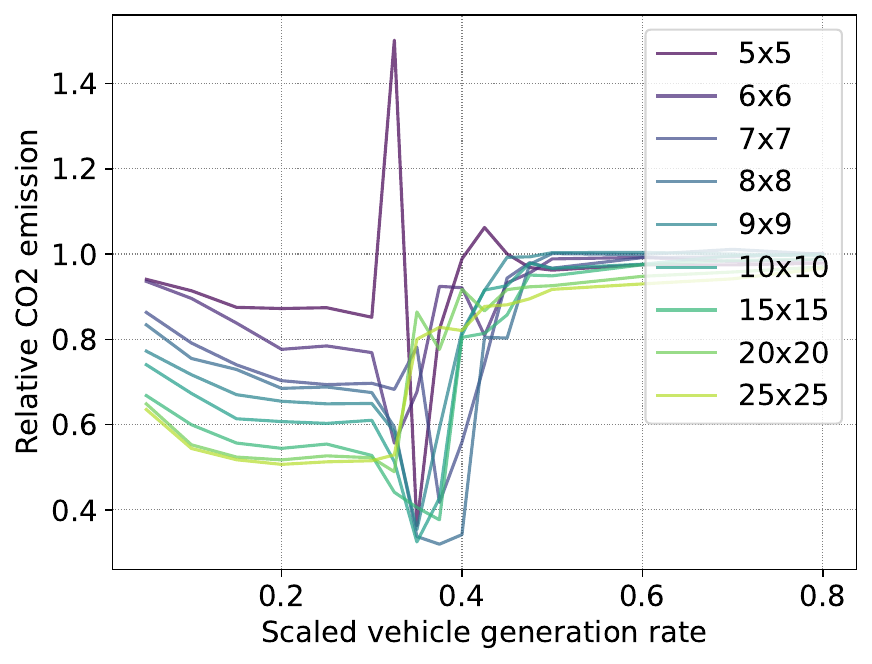}
    \caption{Relative \ce{CO2} emissions.}
    \label{fig:co2-networksize-rel}
  \end{subfigure}
  \caption{Time-averaged performance indicators for various city sizes and vehicle generation rates. Top and bottom panels show absolute and relative values, respectively.
    In the top panel, solid lines show the result of AMPIC and dashed lines show the result of local control. }
  \label{fig:networksize}
\end{figure}

Figure~\ref{fig:networksize} compares the performance of the proposed global control method and local switching control for various vehicle generation rates.
For precision, the horizontal axis is the scaled vehicle generation rate $\tilp$, which is defined as the vehicle generation rate divided by $\sqrt{N}$. This was done because the average distance traveled by an individual vehicle scales as $O(\sqrt{N})$, and thus so does the average time $t_{\rm{av}}$ for which a vehicle remains in the network. In addition, the total road length $L_{\rm{total}}$ of the network scales as $O(N)$. Therefore, the average density of vehicles in the road network is estimated as
$p \times t_{\rm{av}} / L_{\rm{total}} \sim p/\sqrt{N}$, where $p$ is the vehicle generation rate.
This horizontal axis enables us to compare the results for different network sizes. For example, both the mean velocity and the waiting vehicle ratio exhibit universal sudden changes of behavior around $\tilp\sim 0.3$.
Figures~\ref{fig:networksize}d-f plot the results of AMPIC relative to those of the local switching control. Since in the region $\tilp<0.3$ the relative velocity is above unity and the relative waiting vehicle ratio is below unity, AMPIC results in a faster vehicle cruising speed and a lower number of stopped vehicles than those obtained with local switching control. The performance difference increases with network size
and becomes very large
for $\tilp\approx 0.4$ (because congestion occurs only for local switching control).
For high vehicle generation rates ($\tilp> 0.5$), the relative performance difference between the two methods disappears because congestion is unavoidable for both methods.
Similarly to the results for the waiting vehicle ratio, the \ce{CO2} emissions are greatly reduced with the proposed method, especially for a large network.

\subsection{Prediction horizon impact}

\begin{figure}[ht]
  \centering
  \begin{subfigure}{0.32\textwidth}
    \centering\includegraphics[width=0.9\textwidth]{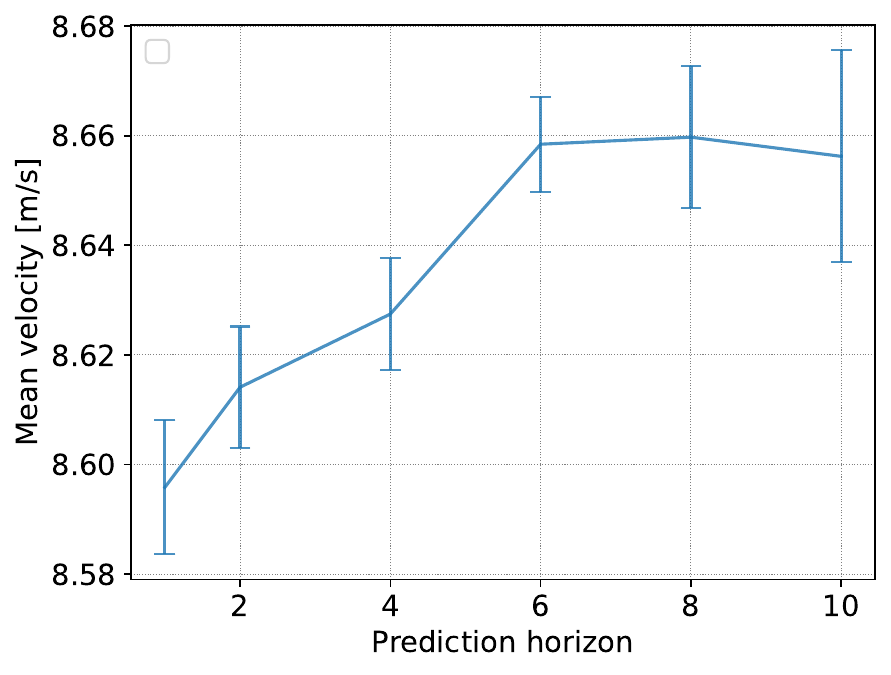}
    \caption{Mean velocity.}
    \label{fig:vmean-horizon}
  \end{subfigure}
  \begin{subfigure}{0.32\textwidth}
    \centering
    \includegraphics[width=0.9\textwidth]{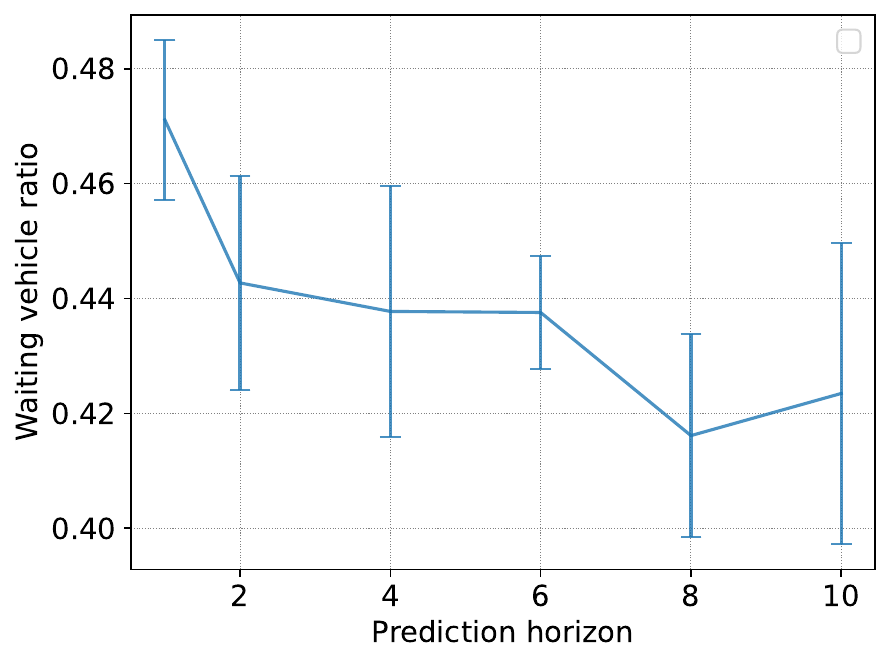}
    \caption{Waiting vehicle ratio.}
    \label{fig:waitrate-horizon}
  \end{subfigure}
  \begin{subfigure}{0.32\textwidth}
    \centering
    \includegraphics[width=0.9\textwidth]{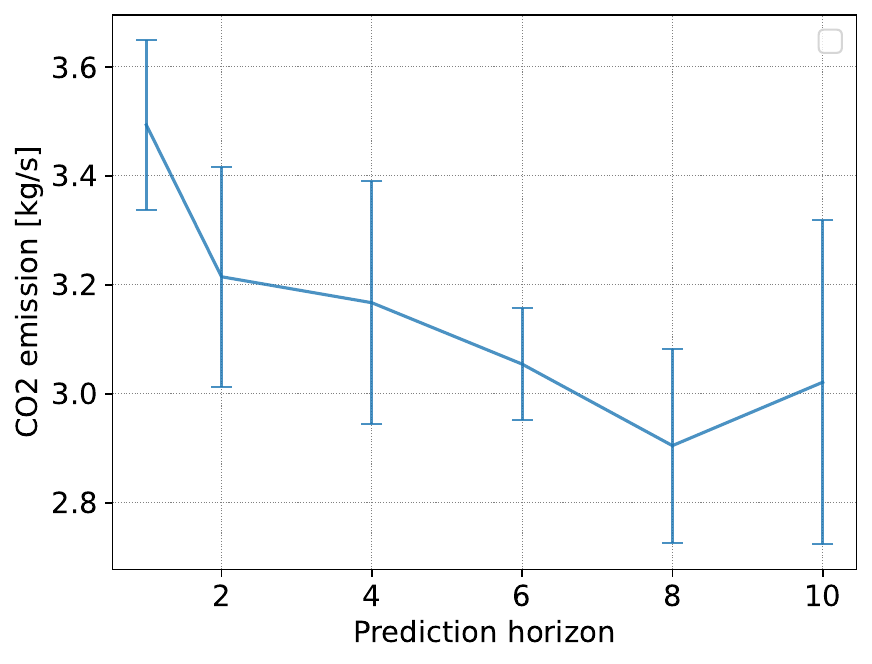}
    \caption{\ce{CO2} emissions.}
    \label{fig:co2-horizon}
  \end{subfigure}
  \caption{Time-averaged performance indicators for various prediction horizons.
    The error bars represent the standard errors of each performance indicator for various seeds of the random number used for vehicle route generation.}
  \label{fig:horizon}
\end{figure}

We next examine the effect of the prediction horizon $k_{\mathrm h}$ on control performance.
Here, we use a lattice network with $N=100$, with the value of $k_{\mathrm h}$ varied from $1$ to $10$.
%We use \emph{SimulatedAnnealingSampler} provided by D-Wave as the Ising solver.
The results are plotted in Fig.~\ref{fig:horizon}.
% The results are plotted in Figs.~\ref{fig:vmean-horizon} and \ref{fig:waitrate-horizon}.
As the prediction horizon is extended, the mean velocity increases and the waiting vehicle ratio tends to decrease, resulting in a decrease in \ce{CO2} emissions.
The slight improvement observed for a horizon longer than $6$ steps indicates the limitation of the model prediction. Given that the control cycle is 60 seconds, this result suggests that considering prediction horizons longer than 6 minutes does not significantly improve traffic conditions.
This limitation may result from the relatively low prediction accuracy for a long horizon, the degradation of optimization performance for a large optimization problem, or both.

\section{Discussion}

We have shown that the model predictive control of traffic signal lights achieved by solving the Ising problem significantly reduces traffic congestion and therefore decreases $\ce{CO2}$ emissions.
The observed improvement compared to the conventional control was higher for a larger traffic network, as shown in Fig.~\ref{fig:networksize}, which is consistent with the fact that all traffic signals in the road network are taken into account simultaneously in the proposed control method. Model predictive control was also shown to work appropriately, as demonstrated for various prediction lengths in Fig.~\ref{fig:horizon}.
Since evaluation was performed using a widely used realistic microscopic traffic flow simulator (SUMO)~\cite{SUMO2018}, implementation in a real city will be possible once some devices are installed to measure the number of vehicles around each intersection to provide important information (i.e., the values of $x_i$ in Eq.~\eqref{eq:def_x}) for the control method.

\begin{table}[ht]
  \centering
  \caption{Time-averaged performance indicators for control implemented with various Ising solvers (Greedy: Greedy method, SA: simulated annealing, QA: quantum annealing).
    % The errors represent the standard errors of each index when the seed of the random number for vehicle route generation is changed.
  }
  \begin{tabular}{lccc}
    \toprule
    {}                                             & {Greedy}                       & {SA}               & {QA}                           \\\midrule
    Mean velocity [\si{\meter\per\second}]         & $8.669 \pm 0.056$              & $8.749 \pm 0.058$  & $8.739 \pm 0.056$              \\
    Waiting vehicle ratio                          & $0.576 \pm 0.025$              & $0.485 \pm 0.024$  & $0.504 \pm 0.025$              \\
    \ce{CO2} emissions [\si{\kilogram\per\second}] & $3.617\pm 0.224$               & $2.722\pm 0.110$   & $2.876 \pm 0.121$              \\
    Sum of squared vehicle bias                    & $93.051 \pm 7.252$             & $56.999 \pm 4.170$ & $68.776 \pm 5.388$             \\
    Elapsed time [\si{\second}]                    & $0.031 \pm 4.64\times 10^{-4}$ & $1.866 \pm 0.168$  & $0.561 \pm 3.42\times 10^{-2}$ \\
    \bottomrule
  \end{tabular}
  \label{table:method}
\end{table}

AMPIC is compatible with all Ising solvers, including those based on quantum annealing. The effect of the solver on performance is of interest
not only for the implementation of traffic signal control but for Ising solver development.
Thus, we compare the results obtained with various Ising solvers for a lattice network with $N=100$.
The following three solvers are used (refer to the supplementary materials for detailed solver settings):
\begin{itemize}
  \item Greedy method: This method is regarded as a discrete analog of the gradient descent method for continuous functions. At each step, the state obtained by flipping the variable that produces the highest energy decay is chosen as the next state.
        % We use \emph{greedy.SteepestDescentSolver} for the specific implementation.
  \item Simulated annealing method (SA): This algorithm searches for the solution by repeatedly transitioning to the next state in a random neighborhood of the current solution.
        The selection is guided by a parameter called the temperature, which is progressively reduced with each iteration of the update.
        This feature reduces the possibility of the solution falling to a local minimum.
        % Its name comes from the process of the same name, in which a metallic material is heated and then gradually cooled to reduce defects in the crystal.
        % We use \emph{SimulatedAnnealingSampler} for the implementation and set the parameter \emph{num\_reads} to $1000$.
  \item Quantum annealing method (QA): This algorithm is regarded as a quantum version of SA. Quantum annealing uses quantum fluctuation to simultaneously search for candidate solutions to find the global minimum of the objective function.
        This is expected to enable faster and more accurate solution searches than those of non-quantum methods.
        % We use \emph{DWaveSampler} for the specific solver and \emph{EmbeddingComposite} for the embeddings for the variable.
        % We set set the parameter \emph{num\_reads} to $1000$.
        % The \emph{Advantage\_system5.4} solver is selected in the D-Wave machine.
\end{itemize}
% 
% The results are shown in Fig.~\ref{fig:method}.
The results are shown in Table \ref{table:method}.
% The results are shown in Figs.~\ref{fig:vmean-method}--\ref{fig:waitrate-method}.
The performance of SA and QA is much better than that of the Greedy method. Particularly, the \ce{CO2} emissions for the Greedy method are more than 25\% higher than those for SA and QA.
The emissions for QA are slightly higher than those for SA, but the
time required to solve the problem with QA is much shorter than that for SA. Note that these results depend on the computational environment.
Using a high-performance CPU for SA will shorten the computation time. Note that \ce{CO2} emitted by the computation itself should also be considered. Such an evaluation is beyond the scope of the present study.
Also note that the computation time for the quantum annealing method was measured as \emph{qpu\_access\_time}~\cite{QPUTime}, which excludes the communication time between Canada (where D-Wave is installed) and Japan (where the experiment was conducted) and the waiting time for the start of processing.
Given the trade-off between performance and time, a quantum annealing machine is a practical candidate for the optimization solver.

The sudden shift from a smooth traffic state to a congested state (see Fig.~\ref{fig:networksize})
is very similar to a phase transition (e.g., that of spin systems).
The large fluctuation (or sample variance shown by the error bars) around the transition point is also compatible with a phase transition.
However, the distinction between the two phases is sharper for small networks, which is
opposite to the general trend for phase-transition phenomena.
This difference could be the result of the intervention of signal control, which acts to blur the phase transition that results in smooth traffic states.
A detailed discussion of this transition phenomenon and an investigation of finite-size scaling for the present system will be presented in future studies.

Improving the overall traffic conditions in large cities is important for achieving carbon neutrality.
In this study, \ce{CO2} emissions were estimated under the assumption that all vehicles are gasoline-fueled. Reducing the power consumption of such vehicles will continue to be crucial even with the increasing prevalence of electric vehicles.
Moreover, the proliferation of connected and automated vehicles is expected to provide more detailed traffic information, improving the accuracy of traffic flow predictions. For the proposed control method, such information can be incorporated in the calculation of parameters $\tilde{A}$ and $\tilde{b}$ in Eq.~\eqref{eq:dynamics_x_vec_in_result}. The inclusion of this information is expected to significantly enhance optimization performance.
% Therefore, this study will contribute to future smart cities with next-generation vehicles.

\section*{Data availability}

The datasets used and analyzed in this study are generated from a code created by the authors.
These datasets can be found on GitHub:
\url{https://github.com/ToyotaCRDL/ampic}.

\section*{Code availability}

The codes for generating all the results can be found on GitHub:
\url{https://github.com/ToyotaCRDL/ampic}.

% \bibliography{my_bib,my_bib_added}

\section*{Acknowledgements}
This study was supported by the Intelligent Mobility Society Design, Social Cooperation Program, a social cooperation research program of Toyota Central R\&D Labs and The University of Tokyo.
The authors would like to thank Dr. Tadayoshi Matsumori and Dr. Norihiro Oyama of Toyota Central R\&D Labs for their useful discussions.

\section*{Author contributions statement}

D.I., H.Ya., and H.Yo. developed the model, carried out the implementation, and analyzed the data.
H.Ya. and H.Yo. developed the computational framework.
D.I. and H.Yo. performed the calculations.
D.I. wrote the manuscript with support from H.Ya., K.A., and H.Yo.
K.A. helped supervise the project.
All authors discussed the results and commented on the manuscript.

\section*{Additional information}

\textbf{Competing interests}: The authors declare no competing financial interests.

\clearpage
\section*{Appendix}

\subsection*{Overview of AMPIC}\label{method}

\begin{algorithm}[ht]
  % \SetAlgoLined
  \caption{Adaptive Model Predictive Ising Controller (AMPIC)}\label{alg:controller}
  \SetKwInOut{Input}{Input}
  \SetKwInOut{Output}{Output}
  \Input{
    $q_{ij}(t)\ (\forall (i,j)\in E)$: Number of vehicles
  }
  \Output{
    $\sigma_{i}(t)\ (\forall i=1,\ldots,N)$: Traffic signal inputs
  }
  Calculate $x(t)$, $a^{\{0,1\}}(t)$, and $o^{\{g,r\}}(t)$ using Eqs.~\eqref{eq-density-normalization}, \eqref{eq:calc_a}, and \eqref{eq-exit-rate}, respectively.\\
  Construct Ising model $C(\sigma([t,\ldots,t+k_h]))$ with Eq.~\eqref{eq:our_ising_model}.\\
  Solve the Ising problem with Ising solver and obtain $\sigma^*([t,\ldots,t+k_h]) = \argmin C(\sigma([t,\ldots,t+k_h]))$.\\
  Return $\sigma^*(t)$.
\end{algorithm}

As shown in Algorithm~\ref{alg:controller}, the proposed controller AMPIC receives the number of vehicles $q(t)$ from SUMO at each time and sends back the traffic signal states in the next control cycle.
In each step, the controller first calculates the necessary parameters for constructing the Ising model.
% Section \ref{sec-normalization}, \ref{sec-exit-rate}.
Using these parameters, the Ising model is constructed.
% , which is described in the following subsection.
% Section \ref{sec-building-qubo}.
Finally, the Ising problem is solved by the Ising solver, and the states of the signals at the current time $t$ are determined.
In the following description of each specific procedure, the function argument $t$ is omitted where it is obvious.

\subsection*{Transformation of the optimal control problem into an Ising problem}\label{sec-building-qubo}

\begin{figure}[ht]
  \centering
  \includegraphics[width=0.8\textwidth]{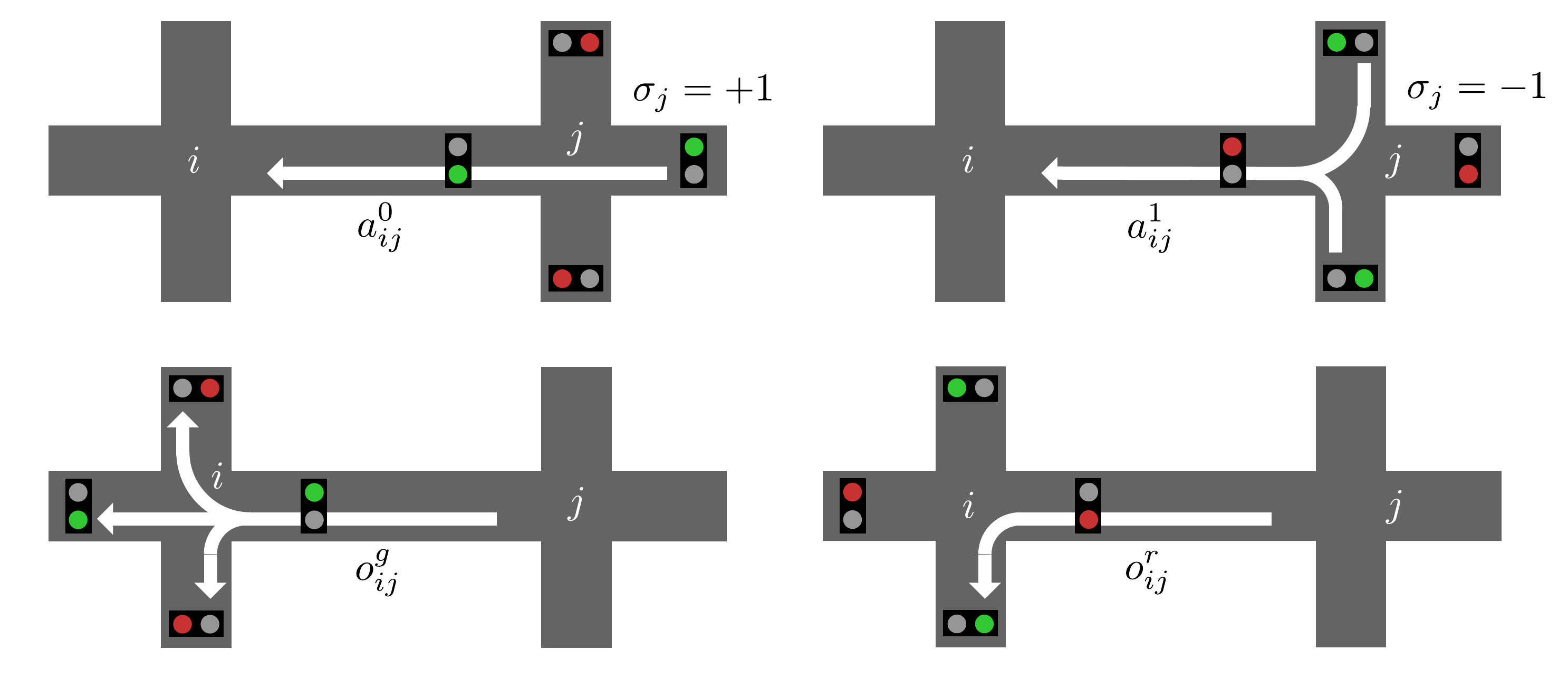}
  \caption{Illustration of each vehicle flow rate on the road network.}
  \label{fig:ao}
\end{figure}

Here, we derive Eq.~\eqref{eq:dynamics_x_vec_in_result} and show that the resulting optimization problem \eqref{eq:eval_func_multi} is represented as an Ising problem.
We first consider the vehicle flow rate at which the cars enter into road $(i,j)$ through intersection $j$ per unit time. Let $a_{ij}^0$ be the flow rate when $\sigma_j=+1$ and $a_{ij}^1$ be that when $\sigma_j=-1$.
We also define $o_{ij}^g$ as the flow rate at which the cars on road $(i,j)$ pass through intersection $i$ per unit time when the traffic light at intersection $i$ is green.
Due to the presence of a dedicated right/left turn lane, vehicles may pass through the intersection $j$ even when the traffic light is red.
In this case, we define $o_{ij}^r$ as the rate of such an outflow.
A visualization of the definition of these flow rates is shown in Fig.~\ref{fig:ao}.
These $a^0_{ij},a^1_{ij},o^g_{ij},o^r_{ij}$ may be known in advance, or they may be adaptively calculated by tracking information on the number of cars $q_{ij}$ while controlling.
The calculation of these values in the present experiment is described later.
%  in Section \ref{sec-exit-rate}.

We can calculate the flow rate at which cars enter and exit road $(i,j)$ as
\begin{align}
   & \frac{1}{2}a^0_{ij}(\sigma_j+1)+\frac{1}{2}a^1_{ij}(-\sigma_j+1)=\frac{1}{2}(\bar{a}_{ij}+a^\Delta_{ij}\sigma_j),                         \\
   & \frac{1}{2}o_{ij}^g (s_{ij}\sigma_i +1)+\frac{1}{2}o_{ij}^r (-s_{ij}\sigma_i +1)=\frac{1}{2}(\bar{o}_{ij}^g+o^\Delta_{ij}s_{ij}\sigma_i),
\end{align}
respectively, where we have defined $(\bar{a}_{ij},a^\Delta_{ij})=(a^0_{ij}+a^1_{ij},a^0_{ij}-a^1_{ij})$ and $(\bar{o}_{ij},o^\Delta_{ij})=(o^g_{ij}+o^r_{ij},o^g_{ij}-o^r_{ij})$.
Then, the time evolution of the number of cars $q_{ij}$ on road $(i,j)$ is
\begin{align}\label{eq:dynamics_q}
  \drv{}{t} q_{ij} =  \frac{1}{2}\bar a_{ij} +\frac{1}{2}a_{ij}^\Delta\sigma_j - \frac{1}{2}\bar{o}_{ij} -\frac{1}{2}o^\Delta_{ij} s_{ij}\sigma_i.
\end{align}
From Eq.~\eqref{eq:dynamics_q} and Eq.~\eqref{eq:def_x}, the time evolution of the vehicle bias $x_i$ is calculated as follows:
\begin{align}
  \begin{split}\label{eq:dynamics_x}
    \drv{}{t} x_i & = \sum_{j\in N_i} \eta_{ij}s_{ij}\qty(\bar a_{ij}+a_{ij}^\Delta\sigma_j  - \bar{o}_{ij} - o^\Delta_{ij} s_{ij}\sigma_i).
  \end{split}
\end{align}
% 
% We rewrite Eq.~\eqref{eq:dynamics_x} in vector form.
Recall that we have defined the vectors $x = [x_1, x_2, \ldots, x_N]^\top$ and ${\sigma} = [\sigma_1, \sigma_2, \ldots, \sigma_N]^\top$.
We can rewrite Eq.~\eqref{eq:dynamics_x} in vector form by using the matrix $A\in\bbR^{N\times N}$ defined as follows:
\begin{align}
  A_{ij} = \begin{cases}
             -\sum_{j\in N_i}\eta_{ij}o_{ij}^\Delta & i=j,                    \\
             \eta_{ij} s_{ij} a_{ij}^\Delta         & i\ne j,\ (i,j)\in E,    \\
             0                                      & i\ne j,\ (i,j)\notin E.
           \end{cases}\label{eq-system-matrix-A}
\end{align}
We also define the vector $ b=[b_1, b_2, \ldots, b_N]^\top\in\bbR^N$ as follows:
\begin{align}
  b_i = \sum_{j\in N_i} \eta_{ij}s_{ij}(\bar a_{ij} - \bar o_{ij}),\quad i=1,\ldots,N.\label{eq-system-matrix-b}
\end{align}
Then, the time evolution in Eq.~\eqref{eq:dynamics_x} is written as follows:
\begin{align}
  \drv{}{t}  x = A{\sigma} + b.
\end{align}
We assumed that the traffic signal states $\sigma(t)$ do not change from time $t$ to $t+\tau$ if the states are updated at time $t$.
Therefore, we can integrate the above to obtain the following difference equation:
\begin{align}
  \begin{split}\label{eq:dynamics_x_vec}
    x(t+\tau) & = x(t) + (A{\sigma}(t) + b)\tau            \\
              & = x(t) + \tilde{A}{\sigma}(t) + \tilde{b},
  \end{split}
\end{align}
where we have defined $\tilde{A}=A\tau$, $\tilde{b}=b\tau$. In the following, we will only consider $t$ on the control cycle: $t=\tau k\ (k\in \bbN)$.

Next, we show that the minimization of the multi-horizon objective function Eq.~\eqref{eq:eval_func_multi} becomes an Ising problem using Eq. ~\eqref{eq:dynamics_x_vec}.
First, we define the vectors $\bfx(t)$ and $\bfsigma(t)$ by using the variables $x(t)$ and $\sigma(t)$ at future $k_h$ points as follows:
\begin{align}
  \bfx(t)     & = [x(t), \ldots, x(t+(k_h-1)\tau)]^\top\in\bbR^{N k_h},           \\
  \bfsigma(t) & = [\sigma(t), \ldots, \sigma(t+(k_h-1)\tau)]^\top\in\bbR^{N k_h}.
\end{align}
We also define the matrices $\bfI, \bfA$ and the vector $\bfb$ as follows:
\begin{align}
  \bfI & = \begin{bmatrix}
             I_N    \\
             I_N    \\
             \vdots \\
             I_N
           \end{bmatrix}\in\bbR^{N k_h\times N},                                 \\
  \bfA & = \begin{bmatrix}
             \tilA  & {}     & {}     & {}              & {}    \\
             \tilA  & \tilA  & {}     & \text{\huge{0}} & {}    \\
             \vdots & \ddots & \ddots & {}              & {}    \\
             \tilA  & \tilA  & \cdots & \tilA           & {}    \\
             \tilA  & \tilA  & \cdots & \tilA           & \tilA \\
           \end{bmatrix}\in\bbR^{N k_h\times N k_h}, \\
  \bfb & = \begin{bmatrix}
             \tilb  \\
             2\tilb \\
             \vdots \\
             k_h \tilb
           \end{bmatrix}\in\bbR^{N k_h},
\end{align}
where $I_N$ denotes the $N$-dimensional identity matrix.
By using Eq.~\eqref{eq:dynamics_x_vec}, we obtain the following equality:
\begin{align}
  \bfx(t+\tau) = \bfI x(t)
  + \bfA \bfsigma(t)
  + \bfb.
\end{align}
Using this, Eq.~\eqref{eq:eval_func_multi} is organized as follows:
\begin{align}
  \begin{split}\label{eq:our_ising_model}
    C(\sigma([t,\ldots,t+(k_h-1)\tau])) & = \sum_{k=1}^{k_h} \qty{  x(t+k\tau)^\top Q x(t+k\tau)}                                                        \\
                                        & = \bfx(t+\tau)^\top\bfQ\bfx(t+\tau)                                                                            \\
                                        & = \qty{\bfI x(t) + \bfA \bfsigma(t) + \bfb }^\top \bfQ\qty{\bfI x(t) +\bfA \bfsigma(t) + \bfb }                \\
                                        & = \bfsigma(t)^\top \bfA^\top \bfQ\bfA \bfsigma(t) + 2 \qty{\bfI x(t)+ \bfb }^\top \bfQ\bfA \bfsigma(t) + \bmc,
  \end{split}
\end{align}
% と整理できる。
where $\bfQ$ is the block diagonal matrix $\bfQ \coloneqq \diag\{Q,\ldots, Q\}$ and $\bmc \coloneqq\qty{\bfI x(t) + \bfb}^\top\bfQ\qty{\bfI x(t) + \bfb}$.
Equation~\eqref{eq:our_ising_model} is nothing but the form of the Ising model of Eq.~\eqref{eq:QUBO}.

\subsection*{Calculation of the weight parameter $\eta$ in the vehicle bias} \label{sec-normalization}

We describe how we determine the coefficients $\eta_{ij}>0$, which was introduced to define vehicle bias in Eq.~\eqref{eq:def_x}.
We consider the reference values for the length and number of cars on the road, denoted by $L^\text{ref}, N^\text{ref}$, respectively.
We then define the normalized density of cars on the road $\tilq_{ij}$ as
\begin{align}
  \tilde{q}_{ij}=\frac{q_{ij}/N^\text{ref}}{L_{ij}/L^\text{ref}},
\end{align}
where $L_{ij}$ denotes the length of the road $(j,i)$.
According to our assumption, there are only four-way or three-way intersections.
When intersection $i$ is four-way, it holds $s_{ij}=+1$ for two of the incoming roads and $s_{ij}=-1$ for the other two.
Therefore, we can easily define the vehicle bias $x_i$ by summing $s_{ij}\tilde{q}_{ij}$ over $j$.

Meanwhile, when intersection $i$ is three-way, either $s_{ij}=+1$ or $s_{ij}=-1$ holds for only one incoming road.
For intersections with such an imbalance, we define the bias of the number of vehicles $x_i$ as follows:
\begin{align}
  x_i = 2\sum_{j\in \mathcal{N}_i} c_{ij} \tilde{q}_{ij},
\end{align}
where the coefficient takes $c_{ij}=2$ for the road $(j,i)$ if there is no other road that leads to the same intersection $i$ and has the same $s_{ij}$, and otherwise $c_{ij}=1$.
From the definition of $\tilde{q}_{ij}$, we obtain
\begin{align}
  x_i = 2\sum_{j\in \mathcal{N}_i} \frac{c_{ij}L^\text{ref}}{N^\text{ref}L_{ij}}q_{ij}, \label{eq-density-normalization}
\end{align}
which leads to the definition of $\eta_{ij}$ as
\begin{align}
  \eta_{ij}=\frac{c_{ij}L^\text{ref}}{N^\text{ref}L_{ij}}.
\end{align}

\subsection*{Calculation of the inflow velocity $a$ and the outflow velocity $o$}\label{sec-exit-rate}

We describe how to determine the rates of inflow $a^{s}_{i,j}$ ($s\in\{0,1\}$) and outflow $o^{s}_{i,j}$ ($s\in\{g,r\}$) for the road $(i,j)$.
First, consider the rate of inflow $a^s_{i,j}$ for $s=0$.
For road $(j,k)$ leading to intersection $j$, when vehicles on the road exit, they will enter one of the roads starting from intersection $j$.
Let $p_{ijk}$ be the rate of vehicles entering such a road $(i,j)$.
Then, the flow rate of vehicles entering road $(i,j)$ through intersection $j$ from road $(j,k)$ can be estimated by multiplying the rate of vehicles exiting road $(j,k)$ by the probability $p_{ijk}$, when intersection $j$ is passable.
Thus, when $s=0$, we obtain
\begin{align}\label{eq:calc_a}
  a^0_{ij}=\sum_{k} o^g_{jk}p_{ijk},
\end{align}
where the sum is taken over $k$ such that the traffic signal in the state $\sigma_{j}=+1$ allows the cars to exit from road $(j,k)$.
In the case of $s=1$, $a^1_{ij}$ is obtained by summing the same values over $k$ such that the outflow is possible when $\sigma_{j}=-1$.
In the experiment in this study, we have all vehicle paths during the simulation are available in advance. Therefore, $p_{ijk}$ can also be calculated before the simulation starts.
The rate of outflow $o^g_{ij}$ can be calculated as $N^{e}_{ij}/T^{g}_{ij}$, where $T^{g}_{ij}$ is the number of seconds when the traffic light indicates green to the road $(i,j)$, and $N^{e}_{ij}$ is the actual number of vehicles that exit the road.
However, $T^g_{ij}$ and $N^{e}_{ij}$ take small values, especially in the transient state of the earlier period of the simulation, and this makes the calculation unstable.
For this reason, in our experiment, the rate of outflow is assumed to be the same for all roads and is calculated as follows:
\begin{align}
  o^g_{ij}=\frac{\sum_{(i, j)\in E}N^e_{ij}}{\sum_{(i, j)\in E} T^g_{ij}}. \label{eq-exit-rate}
\end{align}

Similarly, the rate of outflow $o^r_{ij}$ may depend on the number of vehicles entering the intersection when the traffic light indicates red.
However, it is always $0$ unless there is a dedicated left/right turn lane on the road.

\subsection*{Simulator implementation details}\label{sec:simulator-setting}

All experiments in this study are conducted on a Linux computer with 64 GB of memory and a clock speed of 3.70 GHz.

In SUMO~\cite{SUMO2018}, states such as vehicle positions and speeds are sequentially updated at each second.
To obtain parameters for the Ising problem used by the controller, the simulator collects statistical information about the traffic on the road network.
The road network consists of two lanes, each of which is for opposite directions.
Intersections are assumed to be three- or four-way intersections.
Vehicles are generated at the originating intersection every interval $p$ seconds in simulation time, which means that the vehicle generation rate is $1/p$, and the vehicles are removed from the simulator when they arrive at the destination.
The origin and destination are located at an intersection and are chosen independently at random with uniformly weighted probabilities.
The route plan from origin to destination for each vehicle is generated by \emph{duarouter} (\url{https://sumo.dlr.de/docs/duarouter.html}), a tool included in the SUMO toolset.

In AMPIC, the solvers provided by D-Wave Systems are used.
Specifically, we use \emph{greedy.SteepestDescentSolver} (\url{https://docs.ocean.dwavesys.com/projects/greedy/en/latest/reference/samplers.html}) for the greedy method, \emph{SimulatedAnnealingSampler} (\url{https://docs.ocean.dwavesys.com/projects/neal/en/latest/reference/sampler.html}) for the simulated annealing method, and \emph{DWaveSampler} (\url{https://docs.ocean.dwavesys.com/projects/system/en/stable/reference/samplers.html}) for the quantum annealing method.
For the latter two, the parameter \emph{num\_reads} is set to $1\ 000$.
In the setting of QA, we use \emph{EmbeddingComposite} (\url{https://dwave-systemdocs.readthedocs.io/en/samplers/reference/composites/embedding.html}) for the embeddings for the variables and use the version \emph{Advantage\_system5.4} for the solver.

\end{document}